\documentclass[11pt,a4paper]{article}
\usepackage{amsfonts,amsmath,amssymb,accents,bm,amsthm}%,mathtools}
\usepackage{verbatim}

\usepackage{hyperref}
\usepackage[normalem]{ulem}%for \sout,\xout

\usepackage{bbm}

\usepackage{tikz} %also makes \color work
\usetikzlibrary{math,calc,arrows,positioning,fit,petri,external}

\usepackage{graphicx}
\newcommand{\dis}{\displaystyle}

\newcommand{\noi}{\noindent}
\newcommand{\halmos}{\rule{1ex}{1.4ex}}
\newcommand{\QED}{\nopagebreak{\hspace*{\fill}$\halmos$\medskip}}
\newcommand{\med}{\medskip}
\newcommand{\quand}{\quad\mbox{and}\quad}

\newtheoremstyle{mythm}% name
  {}%      Space above
  {}%      Space below
  {\itshape}%         Body font
  {}%         Indent amount (empty = no indent, \parindent = para indent)
  {\bfseries}% Thm head font
  {}%        Punctuation after thm head
  {.5em}%     Space after thm head: " " = normal interword space;
  {#1 #2 \thmnote{(#3)}}%  Thm head spec (can be left empty, meaning `normal')

\theoremstyle{mythm}
\newtheorem{theorem}{Theorem}%[section]
\newtheorem{proposition}[theorem]{Proposition}
\newtheorem{lemma}[theorem]{Lemma}

\newcommand{\bt}{\begin{theorem}}
\newcommand{\et}{\end{theorem}}
\newcommand{\bl}{\begin{lemma}}
\newcommand{\el}{\end{lemma}}
\newcommand{\bp}{\begin{proposition}}
\newcommand{\ep}{\end{proposition}}
\newenvironment{Proof}[1][]{\noi\textbf{Proof #1}}{\QED}
\newcommand{\bpro}{\begin{Proof}}
\newcommand{\epro}{\end{Proof}}

\newcommand{\be}{\begin{equation}}
\newcommand{\ee}{\end{equation}}
\newcommand{\ba}{\begin{array}}
\newcommand{\ea}{\end{array}}
\newcommand{\bac}{\begin{array}{r@{\,}c@{\,}l}}
\newcommand{\bc}{\be\begin{array}{r@{\,}c@{\,}l}}
\newcommand{\ec}{\end{array}\ee}

\newcommand{\al}{\alpha}

\newcommand{\ga}{\gamma}
\newcommand{\Ga}{\Gamma}
\newcommand{\de}{\delta}
\newcommand{\De}{\Delta}
\newcommand{\eps}{\varepsilon}

\newcommand{\La}{\Lambda}

\newcommand{\tet}{\theta}
\newcommand{\om}{\omega}

\newcommand{\Ci}{{\cal C}}
\newcommand{\Di}{{\cal D}}

\newcommand{\Fi}{{\cal F}}
\newcommand{\Gi}{{\cal G}}
\newcommand{\Hi}{{\cal H}}
\newcommand{\Ii}{{\cal I}}

\newcommand{\Ki}{{\cal K}}
\newcommand{\Li}{{\cal L}}
\newcommand{\Mi}{{\cal M}}
\newcommand{\Ni}{{\cal N}}
\newcommand{\Oi}{{\cal O}}
\newcommand{\Pc}{{\cal P}}

\newcommand{\Ri}{{\cal R}}

\newcommand{\Xc}{{\cal X}}

\newcommand{\E}{{\mathbb E}}

\newcommand{\N}{{\mathbb N}}

\renewcommand{\P}{{\mathbb P}}

\newcommand{\R}{{\mathbb R}}

\newcommand{\Z}{{\mathbb Z}}

\newcommand{\Fb}{{\mathbf F}}

\newcommand{\Rb}{{\mathbf R}}
\newcommand{\Xb}{{\mathbf X}}
\newcommand{\Yb}{{\mathbf Y}}

\newcommand{\rb}{{\mathbf r}}

\newcommand{\desd}{\ensuremath{\Leftrightarrow}}
\newcommand{\volgt}{\ensuremath{\Rightarrow}}
\newcommand{\up}{\uparrow}
\newcommand{\down}{\downarrow}
\newcommand{\sub}{\subset}
\newcommand{\beh}{\backslash}

\newcommand{\asto}[1]{\underset{{#1}\to\infty}{\longrightarrow}}
\newcommand{\Asto}[1]{\underset{{#1}\to\infty}{\Longrightarrow}}

\newcommand{\ti}{\tilde}

\newcommand{\ov}{\overline}
\newcommand{\un}{\underline}

\newcommand{\ffrac}[2]{{\textstyle\frac{{#1}}{{#2}}}}

\newcommand{\di}{\mathrm{d}}
\newcommand{\half}{{[0,\infty)}}

\newcommand{\ha}{\ffrac{1}{2}}

\newcommand{\psim}{\boldsymbol{\psi}_{\rm mon}}

\usepackage{xcolor}

\begin{document}

%numbering formulas within sections
\makeatletter\@addtoreset{equation}{section}
\makeatother\def\theequation{\thesection.\arabic{equation}}

%alternative layout for enumerate lists.
\renewcommand{\labelenumi}{{\rm (\roman{enumi})}}
\renewcommand{\theenumi}{\roman{enumi}}

\title{\vspace*{-2cm}Monotone duality of interacting particle systems}
\author{Jan Niklas Latz\footnote{The Czech Academy of Sciences,
  Institute of Information Theory and Automation,
  Pod vod\'arenskou v\v e\v z\' i~4,
  18200 Praha 8,
  Czech Republic.
  latz@utia.cas.cz, swart@utia.cas.cz.}\textsuperscript{\hspace{0.5em},}\footnote{Charles University, Faculty of Mathematics and Physics}
\and
Jan~M.~Swart${}^\ast$}

\date{\today}

\maketitle

\begin{abstract}\noindent
The duality theory for monotone interacting particle systems was initiated by Gray (1986) and further developed by Sturm and Swart (2018). It contains the better known additive duality as a special case but differs in the sense that the dual process contains not only single particles but also pairs, triples, and general $n$-tuples of particles, which correspond to the fact that in the forward process sometimes several particles are needed to create one particle at a later time. In earlier work, the dual process was constructed for finite initial states only, but, assuming that the empty state is a trap for the forward process, we show that the dual process can be started in infinite initial states and has an upper invariant law. It can therefore be viewed as some sort of interacting particle system in its own right. For the monotone dual of a cooperative contact process, we show that the upper invariant law is the long-time limit started from any nontrivial homogeneous invariant law. We use this to prove continuity of the survival probability of the forward process as a function of its parameters.
\end{abstract}
\vspace{.5cm}

%82C22   	Interacting particle systems in time-dependent statistical
%               mechanics
%60K35   	Interacting random processes; statistical mechanics type models; 
%				percolation theory
{\it Acknowledgements:} Work supported by grant 20-08468S of the Czech Science Foundation (GA\v{C}R). J.N.\ Latz acknowledges additional support from the grant SVV No.\ 260701.

{\setlength{\parskip}{-2pt}\tableofcontents}

\section{Introduction and main results}\label{S:intro}

\subsection{Monotone particle systems}

Let $S$ be a finite set, called the \emph{local state space} and let $\La$ be a countable set, called the \emph{grid}.\footnote{This is often called the \emph{lattice} but we reserve the latter term for its order-theoretic meaning.} Elements $i\in\La$ are called \emph{sites}. Let $S^\La$ denote the space of functions $x:\La\to S$, equipped with the product topology. Elements of $S^\La$ are called \emph{configurations}. A \emph{local map} is a continuous function $m:S^\La\to S^\La$ for which the set
\be
\Di(m):=\big\{i\in\La:\exists x\in S^\La\mbox{ s.t.\ }m(x)(i)\neq x(i)\big\}
\ee
of sites whose state can be changed by $m$ is finite. An \emph{interacting particle system} is a Markov process $(X_t)_{t\geq 0}$ with state space $S^\La$ and generator of the form
\be\label{Gdef}
Gf(x):=\sum_{m\in\Gi}r_m\big\{f\big(m(x)\big)-f\big(x\big)\big\}\qquad(x\in S^\La),
\ee
where $\Gi$ is a countable collection of local maps and $(r_m)_{m\in\Gi}$ are nonnegative rates. Under suitable conditions on the rates (which will be discussed below), it can be shown that the closure of $G$, which is a priori defined for functions depending on finitely many coordinates, is the generator of a Feller process $(X_t)_{t\geq 0}$ on $S^\La$.

To have a concrete example in mind, assume that $S=\{0,1\}$. For each $i,i',j\in\La$, we define local maps by:
\be\ba{l@{\quad}r@{\,}c@{\,}l}
\mbox{(death)}&\dis{\tt dth}_j(x)(k)&:=&\left\{\ba{ll}
0\quad&\mbox{if }k=j,\\
x(k)\quad&\mbox{otherwise,}
\ea\right.\\[20pt]
\mbox{(branching)}&\dis{\tt bra}_{ij}(x)(k)&:=&\left\{\ba{ll}
x(i)\vee x(j)\quad&\mbox{if }k=j,\\
x(k)\quad&\mbox{otherwise,}
\ea\right.\\[20pt]
\mbox{(cooperative branching)}&\dis{\tt coop}_{ii'j}(x)(k)&:=&\left\{\ba{ll}
\big(x(i)\wedge x(i')\big)\vee x(j)\quad&\mbox{if }k=j,\\
x(k)\quad&\mbox{otherwise.}
\ea\right.
\ec
Assume that $(\La,E)$ is a locally finite graph in which each vertex has degree at least two. For each $j\in\La$, we set
\be
\Ni_j:=\big\{i\in\La:\{i,j\}\in E\big\}
\quand
\Ni^{(2)}_j:=\big\{(i,i')\in\La^2:i,i'\in\Ni_j,\ i\neq i'\big\}.
\ee
We will be interested in the interacting particle system with generator
\bc\label{Gcoop}
\dis Gf(x)&:=&\dis(1-\al)\sum_{j\in\La}\frac{1}{|\Ni_j|}\sum_{i\in\Ni_j}\big\{f\big({\tt bra}_{ij}(x)\big)-f\big(x\big)\big\}\\[5pt]
&&\dis+\al\sum_{j\in\La}\frac{1}{|\Ni^{(2)}_j|}\sum_{(i,i')\in\Ni^{(2)}_j}\big\{f\big({\tt coop}_{ii'j}(x)\big)-f\big(x\big)\big\}+\de\sum_{j\in\La}\big\{f\big({\tt dth}_j(x)\big)-f\big(x\big)\big\},
\ec
where $\al\in[0,1]$ and $\de\geq 0$ are model parameters. We call this process the \emph{cooperative contact process} with \emph{cooperation parameter} $\al$ and \emph{death rate} $\de$. For $\al=0$, this is a classical contact process while for $\al=1$ is is purely cooperative.

If the local state space $S$ is equipped with a partial order, then we equip $S^\La$ with the product order and say that a local map $m:S^\La\to S^\La$ is \emph{monotone} if $x\leq y$ implies $m(x)\leq m(y)$. If $S$ is a lattice (i.e., each pair of elements $x,y\in S$ has a least upper bound $x\vee y$ and greatest lower bound $x\wedge y$), then we denote its least element by $0$ and greatest element by $\top$. The configurations that are constantly equal to $0$ or $\top$ are denoted by $\un 0$ and $\un\top$, respectively. We say that a local map $m:S^\La\to S^\La$ is \emph{additive} if $m(\un 0)=\un 0$ and $m(x\vee y)=m(x)\vee m(y)$ $(x,y\in S^\La)$. It is easy to see that additive maps are monotone. In our example, the death and branching maps are additive, while the cooperative branching maps are monotone but not additive. We say that an interacting particle system is monotone\footnote{More precisely, this should be called representably monotone, but we ignore this subtlety.} or additive if its generator can be written in the form (\ref{Gdef}) with each $m\in\Gi$ monotone, or additive, respectively. A \emph{spin system} is an interacting particle system with generator of the form (\ref{Gdef}) such that $S=\{0,1\}$ and $|\Di(m)|=1$ for all $m\in\Gi$.

Duality theory for monotone Markov processes with a totally ordered state space has a long history, see Siegmund \cite{Sie76} and references therein. Interacting particle systems have a state space that is only partially ordered, which necessitates the distinction between additive duality, which is most similar to Siegmund's duality but needs stronger assumptions than monotonicity only, and the more general monotone duality. The duality theory for additive interacting particle systems with local state space $S=\{0,1\}$ was developed by Harris \cite{Har76,Har78} and Griffeath \cite{Gri79}. Later, Gray \cite{Gra86} developed a duality theory for monotone spin systems. Foxall \cite{Fox16} studied additive duality for more general local state spaces while Sturm and Swart \cite{SS18} studied both additive and monotone duality in this more general setting. In the present paper, we further develop the general duality theory for interacting particle systems satisfying the following assumptions:
\be\ba{r@{\ }l}\label{mas}
{\rm(i)}&\mbox{$S$ is a partially ordered set containing a least element $0$.}\\
{\rm(ii)}&\mbox{Each $m\in\Gi$ is monotone with $m(\un 0)=\un 0$.}
\ec
The main novelty of our work is that we allow the dual system to be started in infinite initial states, which allows us to discuss its upper invariant law. To demonstrate the abstract theory, we apply it to the cooperative contact process with generator as in (\ref{Gcoop}). In the remainder of Section~\ref{S:intro} we state our results and discuss them. Our main results are Theorems \ref{T:Fell}, \ref{T:upp}, \ref{T:ergo}, and \ref{T:tetcont}. Proofs will be postponed till Section~\ref{S:proofs}.

\subsection{Graphical representations}

Under suitable conditions on the rates, interacting particle systems with a generator of the form (\ref{Gdef}) can be constructed from a Poisson point process, that in this context is called a graphical representation. We recall this construction here. For each local map $m:S^\La\to S^\La$, we set
\be
\Ri(m):=\big\{(i,j)\in\La^2:\exists x,y\in S^\La\mbox{ s.t.\ }x(k)=y(k)\ \forall k\neq i\mbox{ and }m(x)(j)\neq m(y)(j)\big\},
\ee
and for $i,j\in\La$ we define
\be
\Ri^\up_i(m):=\big\{j\in\La:(i,j)\in\Ri(m)\big\}
\quand
\Ri^\down_j(m):=\big\{i\in\La:(i,j)\in\Ri(m)\big\}.
\ee
Each set $\Ri^\up_i(m)$ is finite since it is contained in $\Di(m)\cup\{i\}$; moreover each $\Ri^\down_j(m)$ is finite by \cite[Lemma~4.13]{Swa22} and the continuity of $m$. Let $\om$ be a Poisson point set on $\Gi\times\R$ with intensity measure $\mu$ defined as
\be
\mu\big(\{m\}\times[s,t]\big):=r_m(t-s)\qquad(m\in\Gi,\ s\leq t).
\ee
We call $\om$ a \emph{graphical representation}. Recall that a function $\half\ni t\mapsto f_t$ is called \emph{cadlag} if it is right continuous with left limits $f_{t-}:=\lim_{s\up t}f_s$ for all $t>0$. As detailed in Appendix~\ref{A:partic}, the following theorem follows from results in \cite{Swa22}. In (\ref{downsum}) below we write $1_A$ for the indicator of a set $A$. In (\ref{fflow}) $1$ denotes the identity map. Then (\ref{fflow}) says that the collection of random maps $(\Xb_{s,t})_{s\leq t}$ form a \emph{stochastic flow}.

\bt[Poisson construction of particle systems]
Assume\label{T:graph} that the rates $(r_m)_{m\in\Gi}$ satisfy
\be\label{downsum}
{\rm(i)}\ \sup_{i\in\La}\sum_{m\in\Gi}r_m1_{\Di(m)}(i)<\infty,\quad
{\rm(ii)}\ \sup_{j\in\La}\sum_{m\in\Gi}r_m\big|\Ri^\down_j(m)\beh\{j\}\big|<\infty.
\ee
Then almost surely, for each $s\in\R$ and $x\in S^\La$, there exists a unique cadlag function $X^{s,x}:[s,\infty)\to S^\La$ such that $X^{s,x}_s=x$ and
\be\label{evo}
X^{s,x}_t=\left\{\ba{ll}
m(X^{s,x}_{t-})\quad&\mbox{if }(m,t)\in\om,\\[5pt]
X^{s,x}_{t-}\quad&\mbox{otherwise}
\ea\right.\qquad(t>s).
\ee
Setting $\Xb_{s,t}(x):=X^{s,x}_t$ $(s\leq t,\ x\in S^\La)$ defines a collection of random continuous maps $(\Xb_{s,t})_{s\leq t}$ from $S^\La$ into itself such that almost surely
\be\label{fflow}
\Xb_{s,s}=1\quand\Xb_{t,u}\circ\Xb_{s,t}=\Xb_{s,u}\quad(s\leq t\leq u).
\ee
If $s\in\R$ and $X_0$ is an $S^\La$-valued random variable, independent of $\om$, then the process $(X_t)_{t\geq 0}$ defined as
\be\label{Xflow}
X_t:=\Xb_{s,s+t}(X_0)\qquad(t\geq 0)
\ee
is a Feller process whose generator is the closure of the operator $G$ defined in (\ref{Gdef}).
\et

We note that the stochastic flow $(\Xb_{s,t})_{s\leq t}$ from (\ref{fflow}) is \emph{stationary} in the sense that $(\Xb_{s+r,t+r})_{s\leq t}$ is equally distributed with $(\Xb_{s,t})_{s\leq t}$, for each $r\in\R$. Using the fact that the restrictions of a Poisson point process to disjoint sets are independent, it is moreover easy to see that $(\Xb_{s,t})_{s\leq t}$ has \emph{independent increments} in the sense that
\be\label{incr}
\Xb_{t_0,t_1},\ldots,\Xb_{t_{n-1},t_n}\quad\mbox{are independent for all }t_0<\cdots<t_n.
\ee
We will need the following result, that will be proved in Appendix~\ref{A:partic}.

\bp[Finite perturbations]
Assume\label{P:perturb} that in addition to (\ref{downsum}),
\be\label{upsum}
\sup_{i\in\La}\sum_{m\in\Gi}r_m\big|\Ri^\up_i(m)\beh\{i\}\big|<\infty.
\ee
Then almost surely, for each $s\leq t$ and $x,y\in S^\La$ such that $\{i\in\La:x(i)\neq y(j)\}$ is finite, the set
\be\label{difset}
\big\{i\in\La:\Xb_{s,t}(x)\neq\Xb_{s,t}(y)\big\}
\ee
is finite.
\ep

\subsection{Duality}

Let $T$ be a finite set and let $\Ci(S^\La,T)$ denote the set of continuous functions $\phi:S^\La\to T$. It turns out that such functions depend on finitely many coordinates \cite[Lemma~4.13]{Swa22} and as a result $\Ci(S^\La,T)$ is countable. Under the assumptions of Theorem~\ref{T:graph}, setting
\be\label{bdef}
\Fb_{t,s}(\phi):=\phi\circ\Xb_{s,t}\qquad\big(t\geq s,\ \phi\in\Ci(S^\La,T)\big)
\ee
defines a collection of random maps $(\Fb_{t,s})_{t\geq s}$ on $\Ci(S^\La,T)$ such that
\be\label{bflow}
\Fb_{s,s}=1\quand\Fb_{t,s}\circ\Fb_{u,t}=\Fb_{u,s}\quad(u\geq t\geq s),
\ee
i.e., $(\Fb_{t,s})_{t\geq s}$ is a \emph{backward stochastic flow}. It is easy to see that $(\Fb_{t,s})_{t\geq s}$ is stationary with independent increments (defined analogously to (\ref{incr})). If $u\in\R$ and $\Phi_0$ is a random variable with values in $\Ci(S^\La,T)$, independent of the graphical representation $\om$, then setting
\be\label{Phidef}
\Phi_t:=\Fb_{u,u-t}(\Phi_0)\quad(t\geq 0)
\ee
defines a continuous-time Markov chain $(\Phi_t)_{t\geq 0}$ with countable state space $\Ci(S^\La,T)$. This continuous-time Markov chain jumps from $\phi$ to $\phi\circ m$ with rate $r_m$, for each $m\in\Gi$. As a consequence of time-reversal one can check that, somewhat unusually, $(\Phi_t)_{t\geq 0}$ has left-continuous sample paths.

From now on, we assume that our interacting particle system satisfies (\ref{mas}). We set
\be\label{xabs}
S^\La_{\rm fin}:=\big\{x\in S^\La:|x|<\infty\big\}
\quad\mbox{with}\quad
|x|:=\big|\{i\in\La:x(i)\neq 0\}\big|\quad(x\in S^\La),
\ee
and we equip $S^\La_{\rm fin}$, which is countable, with the discrete topology. We make the following observation.

\bl[Finite systems]
Assume\label{L:finsys} (\ref{mas}), and assume that the rates $(r_m)_{m\in\Gi}$ satisfy (\ref{downsum}) and (\ref{upsum}). Then, almost surely
\be
\Xb_{s,t}(x)\in S^\La_{\rm fin}\quad\forall s\leq t,\ x\in S^\La_{\rm fin}.
\ee
\el

\bpro
This follows from Proposition~\ref{P:perturb} and the fact that $\Xb_{s,t}(\un 0)=\un 0$ $(s\leq t)$, which is a consequence of (\ref{mas})~(ii).
\epro

We let $\Li(S^\La,\{0,1\})$ denote the space of lower semi-continuous functions $\phi:S^\La\to\{0,1\}$. We say that $\phi:S^\La\to\{0,1\}$ is \emph{monotone} if $x\leq y$ implies $\phi(x)\leq\phi(y)$ and write
\bc
\dis\Li_+(S^\La,\{0,1\})&:=&\dis\big\{\phi\in\Li(S^\La,\{0,1\}):\phi\mbox{ is monotone with }\phi(\un 0)=0\big\},\\[5pt]
\dis\Ci_+(S^\La,\{0,1\})&:=&\dis\big\{\phi\in\Ci(S^\La,\{0,1\}):\phi\mbox{ is monotone with }\phi(\un 0)=0\big\}.
\ec
Combining (\ref{mas})~(ii) and Theorem~\ref{T:graph}, it is easy to prove that the backward stochastic flow $(\Fb_{t,s})_{t\geq s}$ defined in (\ref{bdef}) maps $\Ci_+(S^\La,\{0,1\})$ into itself. Moreover, $(\Fb_{t,s})_{t\geq s}$ can be extended to $\Li_+(S^\La,\{0,1\})$ and maps this space into itself too:

\bl[Preserved subspaces]
Assume\label{L:lowF} (\ref{mas}) and (\ref{downsum}) and for $s\leq t$ and $\phi\in\Li_+(S^\La,\{0,1\})$ define $\Fb_{t,s}(\phi)$ as in (\ref{bdef}). Then, almost surely, for each $s\leq t$ the map $\Fb_{s,t}$ maps the spaces $\Li_+(S^\La,\{0,1\})$ and $\Ci_+(S^\La,\{0,1\})$ into themselves. Moreover, the maps $(\Fb_{t,s})_{t\geq s}$ on $\Li_+(S^\La,\{0,1\})$ form a backward stochastic flow in the sense of (\ref{bflow}).
\el

As in (\ref{Phidef}) we can use the backward stochastic flow $(\Fb_{t,s})_{t\geq s}$ from Lemma~\ref{L:lowF} to define a Markov process $(\Phi_t)_{t\geq 0}$ with state space $\Li_+(S^\La,\{0,1\})$. This Markov process is the main subject of our paper. We start by giving a more concrete description of this process. For any function $\phi:S^\La\to\{0,1\}$, we write
\be
O_\phi:=\big\{x\in S^\La:\phi(x)=1\big\}.
\ee
We note that $\phi\in\Li_+(S^\La,\{0,1\})$ if and only if $O_\phi$ is open, increasing, and not equal to $S^\La$. It will be convenient to characterise $\phi$ by the minimal elements of $O_\phi$. For $x,y\in S^\La$ we write $x<y$ if $x\leq y$ and $x\neq y$. We recall that a \emph{minimal element} of a set $A\subset S^\Lambda$ is a configuration $y\in A$ such that $x\not\in A$ for all $x<y$. For $A\sub S^\La$, we write
\be\label{circ}
A^\circ:=\big\{y\in A:y\mbox{ is a minimal element of }A\big\},\qquad
A^\up:=\big\{x\in S^\La:\exists y\in A\mbox{ s.t.\ }y\leq x\big\}.
\ee
%\bc\label{circ}
%\dis A^\circ&:=&\dis\big\{y\in A:y\mbox{ is a minimal element of }A\big\},\\[5pt]
%\dis A^\up&:=&\dis\big\{x\in S^\La:\exists y\in A\mbox{ s.t.\ }y\leq x\big\}.
%\ec
Then a set $A$ is increasing precisely if it is equal to its ``upset'' $A^\up$. We introduce the spaces
\be\label{Hidef}
\Hi(\La):=\big\{Y\sub S^\La_{\rm fin}:Y^\circ=Y,\ Y\neq\{\un 0\}\big\},\qquad
\Hi_{\rm fin}(\La):=\big\{Y\in\Hi(\La):|Y|<\infty\big\},
\ee
%\bc\label{Hidef}
%\dis\Hi(\La)&:=&\dis\big\{Y\sub S^\La_{\rm fin}:Y^\circ=Y,\ Y\neq\{\un 0\}\big\},\\[5pt]
%\dis\Hi_{\rm fin}(\La)&:=&\dis\big\{Y\in\Hi(\La):|Y|<\infty\big\},
%\ec
where $|Y|$ denotes the cardinality of $Y$. We will prove the following fact. Below $1_{Y^\up}$ is the indicator function of $Y^\up$.

\bp[Encoding monotone lower semi-continuous functions]
The\label{P:openinc} map $Y\mapsto 1_{Y^\up}$ is a bijection from $\Hi(\La)$ to $\Li_+(S^\La,\{0,1\})$ and the map $\phi\mapsto O^\circ_\phi$ is its inverse. Moreover $Y\mapsto 1_{Y^\up}$ is a bijection from $\Hi_{\rm fin}(\La)$ to $\Ci_+(S^\La,\{0,1\})$.
\ep

Combining Lemma~\ref{L:lowF} and Proposition~\ref{P:openinc}, we can define a backward stochastic flow $(\Yb_{t,s})_{t\geq s}$ on $\Hi(\La)$ by
\be\label{Ybdef}
1_{\Yb_{t,s}(Y)^\up}:=\Fb_{t,s}(1_{Y^\up})\qquad\big(t\geq s,\ Y\in\Hi(\La)\big).
\ee
Note that Lemma~\ref{L:lowF} and Proposition~\ref{P:openinc} tell us that this backward stochastic flow maps the space $\Hi_{\rm fin}(\La)$ into itself. Let $\psim:S^\La\times\Hi(\La)\to\{0,1\}$ be defined as
\be\label{psimon}
\psim(x,Y):=1_{Y^\up}(x)\qquad\big(x\in S^\La,\ Y\in\Hi(\La)\big),
\ee
i.e., $\psim(x,Y):=1$ if there exists a $y\in Y$ such that $y\leq x$, and $\psim(x,Y):=0$ otherwise. We claim that the backward stochastic flow $(\Yb_{t,s})_{t\geq s}$ defined in (\ref{Ybdef}) is dual to the stochastic flow $(\Xb_{s,t})_{s\leq t}$ of our interacting particle system with respect to the \emph{duality function} $\psim$ in the sense that
\be\label{dual}
\psim\big(\Xb_{s,t}(x),Y\big)=\psim\big(x,\Yb_{t,s}(Y)\big)\qquad\big(s\leq t,\ x\in S^\La,\ Y\in\Hi(\La)\big).
\ee
Indeed, this follows from the definitions (\ref{Ybdef}) and (\ref{psimon}) by writing
\be
\psim\big(\Xb_{s,t}(x),Y\big)=1_{Y^\up}\circ\Xb_{s,t}(x)=\Fb_{t,s}(1_{Y^\up})(x)=1_{\Yb_{t,s}(Y)^\up}(x)=\psim\big(x,\Yb_{t,s}(Y)\big).
\ee

\subsection{The dual process}

Let $(\Yb_{t,s})_{t\geq s}$ be the backward stochastic flow defined in (\ref{Ybdef}). It is not hard to see that similar to (\ref{Phidef}), $(\Yb_{t,s})_{t\geq s}$ can be used to define a Markov process $(Y_t)_{t\geq 0}$ that we can think of a ``running backwards in time'' compared to the forward process $(X_t)_{t\geq 0}$. In this subsection, we study this Markov process. We start by studying its state space $\Hi(\La)$.

\bp[Dual topology]
There\label{P:topol} exists a unique metrisable topology on $\Hi(\La)$ such that a sequence $Y_n\in\Hi(\La)$ converges to a limit $Y\in\Hi(\La)$ if and only if
\be\label{topol}
1_{Y_n^\up}(x)\asto{n}1_{Y^\up}(x)\quad\forall x\in S^\La_{\rm fin}.
\ee
The space $\Hi(\La)$, equipped with this topology, is compact.
\ep

Our next result says that the Markov process $(Y_t)_{t\geq 0}$  defined by the backward stochastic flow is in fact a Feller process with compact metrisable state space. We call $(Y_t)_{t\geq 0}$ the \emph{monotone dual} of $(X_t)_{t\geq 0}$. Abstract theory tells us that each Feller process is uniquely characterised by its generator. We leave the explicit analytic identification of the generator of $(Y_t)_{t\geq 0}$ for future investigations. In Proposition~\ref{P:backevol} below we will show, however, that $(Y_t)_{t\geq 0}$ can be obtained as the unique solution of an evolution equation similar to (\ref{evo}). We say that a function is \emph{caglad} if it is left continuous with right limits. The fact that the dual process has left-continuous sample paths is a consequence of time reversal.

\bt[Dual process]
Assume\label{T:Fell} (\ref{mas}), and assume that the rates $(r_m)_{m\in\Gi}$ satisfy (\ref{downsum}) and (\ref{upsum}). Let $u\in\R$ and let $Y_0$ be a random variable with values in $\Hi(\La)$, independent of the graphical representation $\om$. Then the process $(Y_t)_{t\geq 0}$ defined as
\be\label{Feldef}
Y_t:=\Yb_{u,u-t}(Y_0)\quad(t\geq 0)
\ee
is a Feller process with caglad sample paths, state space $\Hi(\La)$, and Feller semigroup $(Q_t)_{t\geq 0}$ defined as
\be\label{Qt}
Q_t(Y,\,\cdot\,):=\P\big[\Yb_{0,-t}(Y)\in\,\cdot\,\big]\qquad\big(t\geq 0,\ Y\in\Hi(\La)\big).
\ee
\et

It follows from the remark below (\ref{Ybdef}) that if the dual process $(Y_t)_{t\geq 0}$ is started in an initial state $Y_0\in\Hi_{\rm fin}(\La)$, then $Y_t\in\Hi_{\rm fin}(\La)$ for all $t\geq 0$. In earlier work \cite{Gra86,SS18}, the dual process was only constructed as a continuous-time Markov chain with countable state space $\Hi_{\rm fin}(\La)$. The fact that we allow infinite initial states allows us to define the upper invariant law of $(Y_t)_{t\geq 0}$, which we discuss next.

We equip $\Hi(\La)$ with a partial order such that
\be
\label{dualorder}
Y\leq Z\quad\desd\quad Y^\up\sub Z^\up\qquad(Y,Z\in\Hi(S^\La)).
\ee
The next lemma says that the dual process $(Y_t)_{t\geq 0}$ is monotone with respect to this partial order.

\bl[Monotonicity of the dual process]
Almost\label{L:Ymon} surely, one has $\Yb_{t,s}(Y)\leq\Yb_{t,s}(Z)$ for all $t\geq s$ and $Y,Z\in\Hi(\La)$ such that $Y\leq Z$.
\el

For $a\in S$ and $i\in\La$, we define $e^a_i\in S^\La$ by
\be
e^a_i(j):=\begin{cases}
a&\text{if}\ j=i,\\
0&\text{else},
\end{cases}
\qquad(j\in\Lambda).
\ee
In particular, if $S=\{0,\ldots,n\}$, then we write $e_i:=e^1_i$. We define $Y_{\rm top}\in\Hi(\La)$ by
\be
Y_{\rm top}:=\big\{e^a_i:i\in\La,\ a\in S_{\rm sec}\big\}
\quad\mbox{with}\quad
S_{\rm sec}:=(S\beh\{0\})^\circ.
\ee
In particular, if $S=\{0,\ldots,n\}$, then $S_{\rm sec}=\{1\}$ and $Y_{\rm top}=\{e_i:i\in\La\}$. Elements of $S_{\rm sec}$ are ``second from below'' in the order on $S$, which explains the notation. The next proposition says that $Y_{\rm top}$ is the ``top'' element of $\Hi(\La)$.

\bp[Order on the dual state space]
The\label{P:ord} partial order defined in (\ref{dualorder}) is compatible with the topology on $\Hi(\La)$ in the sense that the set
\be\label{ord}
\big\{(Y,Z)\in\Hi(\La)^2:Y\leq Z\big\}
\ee
is closed with respect to the product topology on $\Hi(\La)^2$. The partially ordered set $\Hi(\La)$ has a least element, which is $\emptyset$, and a greatest element, which is $Y_{\rm top}$.
\ep

As a result of the compatibility condition (\ref{ord}), one can define a stochastic order for probability measures on $\Hi(\La)$ in the usual way, see Appendix~\ref{A:ord}. Together with Lemma~\ref{L:Ymon} this allows us to apply an abstract result (Proposition~\ref{P:upp} in the appendix) to draw the following conclusions.

\bt[Upper invariant law]
Assume\label{T:upp} (\ref{mas}), and assume that the rates $(r_m)_{m\in\Gi}$ satisfy (\ref{downsum}) and (\ref{upsum}). Then the Feller process $(Y_t)_{t\geq 0}$ with semigroup $(Q_t)_{t\geq 0}$ defined in (\ref{Qt}) has an invariant law $\ov\nu_{\rm Y}$ with the property that if $\nu$ is another invariant law, then $\nu\leq\ov\nu_{\rm Y}$ in the stochastic order. Moreover, the process started in $Y_{\rm top}$ satisfies
\be\label{Yup}
\P^{Y_{\rm top}}\big[Y_t\in\,\cdot\,\big]\Asto{t}\ov\nu_{\rm Y},
\ee
where $\Rightarrow$ denotes weak convergence of probability measures on $\Hi(\La)$.
\et

More trivially, the process $(Y_t)_{t\geq 0}$ also has a lower invariant law, but since this is simply the delta measure concentrated on the least element of $\Hi(\La)$, which is $\emptyset$, this lower invariant law is less interesting. It is easy to see that the partial order and the topology on $S^\La$ satisfy a condition similar to (\ref{ord}) and hence are also ``compatible''. As a result, if the local state space $S$ has a greatest element $\top$, then we can apply Proposition~\ref{P:upp} in the appendix to conclude that under the assumptions of Theorem~\ref{T:upp}, the interacting particle system $(X_t)_{t\geq 0}$ also has an upper invariant law $\ov\nu_{\rm X}$, which is given by
\be\label{Xup}
\P^{\un\top}\big[X_t\in\,\cdot\,\big]\Asto{t}\ov\nu_{\rm X},
\ee
where $\un\top$ denotes the configuration that is identically equal to $\top$ and $\Rightarrow$ denotes weak convergence of probability measures on $S^\La$. This upper invariant law dominates all other invariant laws of $(X_t)_{t\geq 0}$ in the stochastic order.

\subsection{Survival and stability}\label{S:surv}

Throughout this subsection we assume (\ref{mas}), (\ref{downsum}), and (\ref{upsum}). We also assume that $S$ has a greatest element $\top$ so that the upper invariant law $\ov\nu_{\rm X}$ of the forward process is well-defined. We say that the interacting particle system $(X_t)_{t\geq 0}$ is \emph{stable} if $\ov\nu_{\rm X}\neq\de_{\un 0}$, the delta measure on the all zero configuration. It is easy to see that this implies that $\ov\nu_{\rm X}(\{\un 0\})=0$. Indeed, if this would not be the case then we could write $\ov\nu_{\rm X}=p\nu'+(1-p)\de_{\un 0}$ for some $0<p<1$, where $\nu'$ would be an invariant law that is larger than $\ov\nu_{\rm X}$ in the stochastic order, contradicting the maximality of the latter. We say that the dual process $(Y_t)_{t\geq 0}$ is \emph{stable} if $\ov\nu_{\rm Y}\neq\de_\emptyset$, the delta measure on $\emptyset$, which is the least element of $\Hi(\La)$. By the same argument as before, this implies that $\ov\nu_{\rm Y}(\{\emptyset\})=0$.

We say that the interacting particle system $(X_t)_{t\geq 0}$, respectively its monotone dual $(Y_t)_{t\geq 0}$ \emph{dies out} if
\be\ba{ll}
\dis\P\big[\exists t\geq 0\mbox{ s.t.\ }\Xb_{0,t}(x)=\un 0\big]=1&\dis\quad(x\in S^\La_{\rm fin}),\\[5pt]
\dis\P\big[\exists t\geq 0\mbox{ s.t.\ }\Yb_{0,-t}(Y)=\emptyset\big]=1&\dis\quad\big(Y\in\Hi_{\rm fin}(\La)\big).
\ec
If these probabilities are less than one for some $x\in S^\La_{\rm fin}$ or $Y\in\Hi_{\rm fin}(\La)$, then we say that $(X_t)_{t\geq 0}$ or $(Y_t)_{t\geq 0}$ \emph{survives}. The following lemma links these concepts to stability.

\bl[Survival and stability]
Assume\label{L:surstab} (\ref{mas}), (\ref{downsum}), and (\ref{upsum}), and assume that $S$ has a greatest element $\top$. Let $(X_t)_{t\geq 0}$ be the interacting particle system with generator (\ref{Gdef}) and let $(Y_t)_{t\geq 0}$ be its monotone dual. Then:
\begin{itemize}
\item[{\rm(i)}] $(X_t)_{t\geq 0}$ is stable if and only if $(Y_t)_{t\geq 0}$ survives.
\item[{\rm(ii)}] $(X_t)_{t\geq 0}$ survives if and only if $(Y_t)_{t\geq 0}$ is stable.
\end{itemize}
\el

We note that part~(i) of the lemma (specialised to spin systems) already occurs as formula (30) in \cite{Gra86}. Part~(ii) had to wait to the present paper, since $(Y_t)_{t\geq 0}$ has (to the best of our knowledge) not previously been constructed for infinite initial states.

%Gray 86 Section 2 does stochastic flows (notation \xi(s,t,A)).
%Dual stochastic flow denoted \zeta.
%Gray does not have the condition m(\un 0)=\un 0 so his formula (29)
%is more general but in (30) he specialises to the case that m(\un 0)=\un 0.

This is a good moment to demonstrate the abstract theory developed so far on the cooperative contact process with generator as in (\ref{Gcoop}). For concreteness, we look at the system on $\La=\Z^d$ with nearest neighbour edges. We let
\be
\rho(\al,\de):=\int\!\ov\nu_{\rm X}(\di x)\,x(i)
\quand
\tet(\al,\de):=\P^{e_i}\big[X_t\neq\un 0\ \forall t\geq 0\big]
\qquad(i\in\Z^d)
\ee
denote the density of the upper invariant law and the survival probability started from a configuration containing a single one, respectively. Note that by translation invariance, these quantities do not depend on $i\in\Z^d$. Figure~\ref{fig:coop} shows the conjectured phase diagram in dimension $d=2$. The phase diagram is believed to be similar in higher dimensions, but not in one dimension. Figure~\ref{fig:coop} is based on numerical data for the quantities $\rho(\al,\de)$ and $\tet(\al,\de)$, partial rigorous results for this and a related model, and heuristics based on the mean-field equation. Simple coupling arguments (see Lemma~\ref{L:monnu} below) show that $\rho$ and $\tet$ are nonincreasing in $\al$ and $\de$. As a result, there exist $\de_{\rm c}(\al)\geq 0$ and $\de'_{\rm c}(\al)\geq 0$ that are nonincreasing as a function of $\al\in[0,1]$ such that
\be\label{dec}
\rho(\al,\de)\left\{\ba{ll}
>0\quad&\mbox{for }\de<\de_{\rm c}(\al),\\
=0\quad&\mbox{for }\de>\de_{\rm c}(\al),
\ea\right.
\quand
\tet(\al,\de)\left\{\ba{ll}
>0\quad&\mbox{for }\de<\de'_{\rm c}(\al),\\
=0\quad&\mbox{for }\de>\de'_{\rm c}(\al).
\ea\right.
\ee
Numerical data suggest the existence of a critical value $0<\al_{\rm c}<1$ such that $\de'_{\rm c}(\al)=\de_{\rm c}(\al)$ for $\al\leq\al_{\rm c}$ and $\de'_{\rm c}(\al)<\de_{\rm c}(\al)$ for $\al>\al_{\rm c}$.

\begin{figure}[htb]
\begin{center}
\begin{tikzpicture}[>=triangle 45]
\draw[thick] (0,5)--(7.6,1.35);
\draw[thick] (7.6,1.35)--(9,0.75);
\draw[thick] (7.6,1.35)--(9,0);
\draw (4.2,1.5) node[black] {\parbox{7cm}{survival and stability}};
\draw (6,5) node[black] {\parbox{7cm}{extinction and instability}};
\draw (7.5,2.8) node[black] {\parbox{2.5cm}{extinction\\ and stability}};
\draw[black,very thick,->] (7.5,2.3)--(8.7,0.6);

\draw (4,3.7) node {$\de'_{\rm c}=\de_{\rm c}$};
\draw (8.1,0.4) node {$\de'_{\rm c}$};
\draw (8.7,1.1) node {$\de_{\rm c}$};

\draw (0,0)--(9,0)--(9,6)--(0,6)--(0,0);
\draw[very thick,->] (3,-0.7)--+(2,0) node[right]{$\al$};
\draw[very thick,->] (-0.9,2)--+(0,2) node[above]{$\de$};
\draw[very thick] (0,1.71)--+(-0.1,0) node[left]{$0.2$};
\draw[very thick] (0,3.43)--+(-0.1,0) node[left]{$0.4$};
\draw[very thick] (0,5.14)--+(-0.1,0) node[left]{$0.6$};
\draw[very thick] (1.8,0)--+(0,-0.1) node[below]{$0.2$};
\draw[very thick] (3.6,0)--+(0,-0.1) node[below]{$0.4$};
\draw[very thick] (5.4,0)--+(0,-0.1) node[below]{$0.6$};
\draw[very thick] (7.2,0)--+(0,-0.1) node[below]{$0.8$};
\draw[very thick] (9,0)--+(0,-0.1) node[below]{$1$};
\end{tikzpicture}
\caption{Conjectured phase diagram for the cooperative contact process on $\Z^2$. The density of the upper invariant law $\rho(\al,\de)$ is positive for $\de<\de_{\rm c}(\al)$ and the survival probability $\tet(\al,\de)$ is positive for $\de<\de'_{\rm c}(\al)$. Numerically, one sees that $\rho(\al,\de)>0$ also for $\de=\de_{\rm c}(\al)$ in the regime where $\de'_{\rm c}(\al)<\de_{\rm c}(\al)$ while the functions $\rho$ and $\tet$ are continuous everywhere else.}
\label{fig:coop}
\end{center}
\end{figure}

Let $(\Xb_{s,t})_{s\leq t}$ be the stochastic flow of the cooperative contact process $(X_t)_{t\geq 0}$ and let $(\Yb_{t,s})_{t\geq s}$ be the backward stochastic flow associated with its monotone dual $(Y_t)_{t\geq 0}$. Recall that the dual process $(Y_t)_{t\geq 0}$ takes values in the space $\Hi(\La)$ of subsets $Y\sub S^\La_{\rm fin}$ that satisfy $Y^\circ=Y$ and $Y\neq\{\un 0\}$. The duality relation (\ref{dual}) says that
\be\label{du1}
\exists y\in Y\mbox{ s.t.\ }y\leq\Xb_{s,t}(x)
\quad\desd\quad
\exists y\in\Yb_{t,s}(Y)\mbox{ s.t.\ }y\leq x
\ee
for all $s\leq t$, $x\in S^\La$, and $Y\in\Hi(\La)$. Taking into account that $\Yb_{t,s}(Y)^\circ=\Yb_{t,s}(Y)$, this mean that:
\be\ba{c}
\dis\Yb_{t,s}(Y)\mbox{ is the set of minimal elements of }\big\{x\in S^\La_{\rm fin}:\Xb_{s,t}(x)\in A\big\}\\[5pt]
\mbox{where }A:=\big\{x'\in S^\La_{\rm fin}:\exists y\in Y\mbox{ s.t.\ }y\leq x'\big\}.
\ec
In particular, setting $Y=\{e_i\}$ we have $A=\{x'\in S^\La_{\rm fin}:x'(i)=1\}$ and $\Yb_{t,s}(Y)$ is the set of minimal configurations required at time $s$ to create a particle at $i$ a time $t$. Letting $\un 1$ denote the configuration that is identically one, one has, for any $i\in\La$,
\be
\P^{\un 1}\big[X_t(i)=1\big]=\E\big[\psim\big(\Xb_{0,t}(\un 1),e_i\big)\big]=\E\big[\psim\big(\un 1,\Yb_{t,0}(e_i)\big)\big]=\P^{e_i}\big[Y_t\neq\emptyset\big],
\ee
so taking the limit $t\to\infty$ we see that
\be
\rho(\al,\de)=\P^{e_i}\big[Y_t\neq\emptyset\ \forall t\geq 0\big],
\ee
proving Lemma~\ref{L:surstab}~(i) for the cooperative contact process. The proof of Lemma~\ref{L:surstab}~(ii) is similar. For any $i\in\La$, one has
\be
\P^{e_i}\big[X_t\neq\un 0\big]=\E\big[\psim\big(\Xb_{0,t}(e_i),Y_{\rm top}\big)\big]=\E\big[\psim\big(e_i,\Yb_{t,s}(Y_{\rm top})\big)\big]=\P^{Y_{\rm top}}\big[e_i\in Y_t\big],
\ee
which after taking the limit $t\to\infty$ yields
\be\label{tetovnu}
\tet(\al,\de)=\P^{e_i}\big[X_t\neq\un 0\ \forall t\geq 0\big]=\int_{\Hi(\La)}\ov\nu_{\rm Y}(\di Y)1_{\{e_i\in Y\}}.
\ee

We can give a more explicit description of the evolution of the process $(Y_t)_{t\geq 0}$ as follows. In general, if $S$ is a partially ordered set containing a least element $0$ and $m:S^\La\to S^\La$ is a continuous map satisfying $m(\un 0)=\un 0$, then we will show in Lemma~\ref{L:dualmap} below that there exists a unique dual map $\hat m:\Hi(\La)\to\Hi(\La)$ such that
\be
\psim\big(m(x),Y\big)=\psim\big(x,\hat m(Y)\big)\qquad\big(s\leq t,\ x\in S^\La,\ Y\in\Hi(\La)\big).
\ee
The formal generator of $(Y_t)_{t\geq 0}$ is then given by
\be
\hat Gf(Y):=\sum_{m\in\Gi}r_m\big\{f\big(\hat m(Y)\big)-f\big(Y\big)\big\}\qquad\big(Y\in\Hi(\La)\big).
\ee
In the special case of the cooperative contact process, one can check that
\be\ba{c}\label{dumaps}
\dis\widehat{\tt dth}_j(Y)=\{y\in Y:y(j)=0\},\quad
\widehat{\tt bra}_{ij}(Y)=\big(Y\cup\{(y-e_j)\vee e_i:y\in Y,\ y(j)=1\}\big)^\circ\\[5pt]
\widehat{\tt coop}_{ii'j}(Y)=\big(Y\cup\{(y-e_j)\vee e_i\vee e_{i'}:y\in Y,\ y(j)=1\}\big)^\circ.
\ec
We observe that the maps $\widehat{\tt dth}_j$ and $\widehat{\tt bra}_{ij}$ have the property that they map the space
\be
\Hi_{\rm add}(\La):=\big\{Y\in\Hi(\La):|y|=1\ \forall y\in Y\big\}
\ee
into itself. As a consequence, for $\al=0$ the monotone dual $(Y_t)_{t\geq 0}$ of the cooperative contact process has the property that $Y_0\in\Hi_{\rm add}(\La)$ implies $Y_t\in\Hi_{\rm add}(\La)$ for all $t\geq 0$. We claim that this is a consequence of the fact that the process with $\al=0$ is additive. Indeed, we can naturally identify $\{0,1\}^\La$ with $\Hi_{\rm add}(\La)$ via the bijection
\be
\{0,1\}^{\Z^d}\ni z\mapsto\big\{e_i:z(i)=1\big\}\in\Hi_{\rm add}(\La).
\ee
Identifying $Y_t$ with an element of $\{0,1\}^\La$ in this way, one can check that $(Y_t)_{t\geq 0}$ is the additive dual of $(X_t)_{t\geq 0}$ in the sense of additive systems duality. In particular, because of the well-known self-duality of the contact process, $(Y_t)_{t\geq 0}$ is in fact a contact process. Thus, one can view the monotone duals considered in this article as an extension of the classical additive duals of additive interacting particle systems.

\subsection{Ergodicity of the dual process}

We continue to look at the cooperative contact process on $\Z^d$ with nearest neighbour edges. We define translation operators $S_i:\{0,1\}^{\Z^d}\to\{0,1\}^{\Z^d}$ and $S_i:\Hi(\Z^d)\to\Hi(\Z^d)$ by
\be
(S_ix)(j):=x(j-i)\quad(j\in\Z^d)\quand
(S_iY):=\big\{S_iy:y\in Y\big\}.
\ee
We say that a probability law $\mu$ on $\Hi(\Z^d)$ is \emph{homogeneous} if it is invariant under translations, i.e., $\mu\circ S^{-1}_i=\mu$ for all $i\in\Z^d$. We will prove the following result.

\bt[Homogeneous initial laws]
Assume\label{T:ergo} that $(\al,\de)\neq(1,0)$. Let $(Y_t)_{t\geq 0}$ be the monotone dual of the cooperative contact process on $\Z^d$ started in an initial law that is homogeneous and satisfies $Y_0\neq\emptyset$ a.s. Then
\be\label{ergo}
\P\big[Y_t\in\,\cdot\,\big]\Asto{t}\ov\nu_{\rm Y},
\ee
where $\Rightarrow$ denotes weak convergence of probability measures on $\Hi(\Z^d)$. 
\et

As a consequence of Theorem~\ref{T:ergo}, by a well-known argument based on the relation (\ref{tetovnu}), we can prove the following fact about the process we were originally interested in.

\bt[Continuity of the survival probability]
The\label{T:tetcont} function $(\al,\de)\mapsto\tet(\al,\de)$ is upper semi-continuous on $[0,1]\times\half$, and continuous on the set
\be\label{alde}
A:=\big\{(\al,\de)\in[0,1)\times\half:\exists\eps>0\mbox{ s.t.\ }\tet(\al+\eps,\de+\eps)>0\big\}.
\ee
\et

\subsection{Discussion and open problems}

We have extended the duality theory for monotone interacting particle systems to allow the dual process to be started in infinite initial states, and used this to prove a property of one particular monotone interacting particle system, the cooperative contact process. Apart from monotonicity, we also assumed that the all-zero configuration $\un 0$ is a trap for the interacting particle system. As long as one is only interested in the dual process started in finite initial states, this condition can be dropped, as already shown by Gray \cite{Gra86}. In this case, the state space $\Hi_{\rm fin}(\La)$ needs to be extended so that it also contains $\{\un 0\}$ (compare (\ref{Hidef})), which now becomes a trap for the dual process. In principle, this construction also works for infinite initial states but it seems that in most cases of interest, the dual process started in an infinite initial state now jumps to the trap $\{\un 0\}$ immediately. Only when $\un 0$ is a trap for the interacting particle system, the trap $\{\un 0\}$ for the dual process becomes inaccessible, which is why we have made this assumption throughout.

There are plenty of open problems for monotone interacting particle systems that are not additive, and for the cooperative contact process in particular. Here we only briefly discuss the latter. For $\al=0$ the process reduces to a standard contact process, which is self-dual. As a consequence, $\rho(0,\de)=\tet(0,\de)$ for all $\de\geq 0$ and therefore $\de'_{\rm c}(0)=\de_{\rm c}(0)$. It is well-known that $0<\de_{\rm c}(0)<\infty$ \cite{Lig99}. Since $\rho$ and $\tet$ are nonincreasing in $\al$ it follows that $\de_{\rm c}(\al)$ and $\de'_{\rm c}(\al)$ are finite for all $\al\in[0,1]$. A simple monotone coupling with a process without cooperative branching, together with a rescaling of time, tells us that $\de_{\rm c}(\al)\geq(1-\al)\de_{\rm c}(0)$ and likewise for $\de'_{\rm c}(\al)$, so both functions are $>0$ for all $\al<1$. It is easy to see that $\de'_{\rm c}(1)=0$, since the process with $\al=1$, started in a finite interval, rectangle, or (hyper)cube cannot escape from such a set and hence a.s.\ dies out for each $\de>0$. In dimensions $d\geq 2$ it is known that $\de_{\rm c}(1)>0$. This follows from \cite[Thm~18.3.1]{Gra99}. The proof is quite hard and does not yield a lower bound on $\de_{\rm c}(1)$ that is anywhere near the numerically observed value.

Well-known results for the contact process \cite{Lig99}  tell us that the function $\de\mapsto\tet(0,\de)$ is continuous, with the most difficult statement, continuity at $\de_{\rm c}(0)$, having been proved in \cite{BG90}. The methods of the latter paper were generalised in \cite{BG94}. Applying \cite[Thm~2.4]{BG94} to the cooperative contact process one finds that:
\be
\mbox{The set }\big\{(\al,\de)\in[0,1]\times(0,\infty):\tet(\al,\de)>0\big\}\mbox{ is an open subset of }[0,1]\times(0,\infty),
\ee
which combined with our Theorem~\ref{T:tetcont} shows that $\tet$ is continuous everywhere except in the point $(\al,\de)=(1,0)$. The methods of \cite{BG94} can also be used to show that $\tet(\al,\de)>0$ implies $\rho(\al,\de)>0$ and hence $\de'_{\rm c}(\al)\leq\de_{\rm c}(\al)$ for all $\al\in[0,1]$. This is folklore; the details have unfortunately never been written down. See \cite[Thm~1.2]{Fox15}, however, for a sketch of the argument for a different model.

Apart from the facts we have just mentioned, little seems to be known. It is not known if $\de'_{\rm c}(\al)=\de_{\rm c}(\al)$ for $\al$ sufficiently small, as suggested by numerical simulations of the process on $\Z^2$, nor do there seem to be known results that would allow us to conclude that $\de'_{\rm c}(\al)<\de_{\rm c}(\al)$ for $\al$ sufficiently close to one. The continuity of $\rho$ on the set $\{(\al,\de):\rho(\al,\de)>0\}$ also seems to be an open problem, except on the strip $\al=0$. This continuity would follow if one could prove that all homogeneous invariant laws are convex combinations of $\ov\nu_{\rm X}$ and the delta measure on $\un 0$. To prove this, it would suffice to prove in analogy with Theorem~\ref{T:ergo} that the process $(X_t)_{t\geq 0}$ started in arbitrary nontrivial homogeneous initial law $\mu$ converges to $\ov\nu_{\rm X}$. The latter statement is well-known for the contact process \cite{Har76} but it seems plausible it is in general \emph{not true} for the cooperative contact process. Let $\P^{\pi_p}$ denote the law of the process started in product law $\pi_p$ with intensity $p$, and let
\be
\eta(\al,\de):=\inf\big\{p\in[0,1]:\P^{\pi_p}\big[X_t\in\,\cdot\,\big]\Asto{t}\ov\nu_{\rm X}\big\}
\ee
denote the minimal initial density for which the process started in product law converges to the upper invariant law. Then we conjecture that $\eta(\al,\de)>0$ in the regime $\de'_{\rm c}(\al)<\de\leq\de_{\rm c}(\al)$. Proving this is an open problem. Finally, numerical simulations suggest that $\de\mapsto\rho(\al,\de)$ is continuous at $\de_{\rm c}(\al)$ if $\de'_{\rm c}(\al)=\de_{\rm c}(\al)$, but not if $\de'_{\rm c}(\al)<\de_{\rm c}(\al)$. Proving this, except at $\al=0$, is an open problem too. As a last remark, it is worth mentioning that \cite{BG94} and \cite{Gra99} use the structure of $\Z^d$ in an essential way and cannot easily be generalised to different grids, so that on more general grids even less is known.

At this moment, it is too early to tell if monotone duality can help resolve some of these open problems. An interesting line of thought seems to be to try and prove that for $\al$ sufficiently small, the dual process behaves almost as in the additive case $\al=0$, for example in the sense that if $Y$ is distributed according to the upper invariant law, then elements $y\in Y$ with $|y|\geq 2$ have a low density. Another question is under what assumptions the homogeneity assumption in Theorem~\ref{T:ergo} can be dropped, i.e., if one can prove some form of \emph{complete convergence} for the dual process. (See \cite{Lig99} for this terminology in the context of the contact process.)

\subsubsection*{Outline}

The remainder of the paper is devoted to proofs. Propositions \ref{P:openinc}, \ref{P:topol}, and \ref{P:ord} are proved in Subsection~\ref{S:H(S)}. Lemma~\ref{L:lowF} and Theorem~\ref{T:Fell} are proved in Subsection~\ref{S:DualMP}. Lemma~\ref{L:Ymon}, Theorem~\ref{T:upp}, and Lemma~\ref{L:surstab} are proved in Subsection~\ref{S:UP&Surv}. In Subsection~\ref{S:homlaw} we prove a more general version of Theorem~\ref{T:ergo} that is then in Subsection~\ref{S:coop} used to prove Theorems \ref{T:ergo} and \ref{T:tetcont}. The proofs of Theorem~\ref{T:graph} and Proposition~\ref{P:perturb} can be found in Appendix~\ref{A:partic}.

\section{Proofs}\label{S:proofs}

\subsection{The dual space}\label{S:H(S)}

In this section we prove Propositions \ref{P:openinc}, \ref{P:topol}, and \ref{P:ord}.\med

\bpro[Proof of Proposition~\ref{P:openinc}]
We claim that for any $A\sub S^\La$,
\be\label{elem}
{\rm(i)}\ (A^\circ)^\circ=A^\circ,\quad
{\rm(ii)}\ (A^\up)^\up=A^\up,\quad
{\rm(iii)}\ (A^\up)^\circ=A^\circ.
\ee
Indeed, the first two properties are trivial while (iii) follows from (i) and the observations that $A^\circ\sub(A^\up)^\circ\sub A$. We next claim that for $A\sub S^\La$
\be\label{cirup}
A^\up=A,\ A\mbox{ open}\quad\volgt\quad A^\circ\sub S^\La_{\rm fin}\mbox{ and }(A^\circ)^\up=A.
\ee
To see this, we first note that trivially $(A^\circ)^\up\sub A^\up=A$. To prove the converse inclusion, assume that $x\in A$. Then we can find $x_n\in S^\La_{\rm fin}$ such that $x_n\to x$. Since $A$ is open there exists an $n$ such that $x_n\in A$. Since $x_n\in S^\La_{\rm fin}$ we can find, in a finite number of steps, an $x'\leq x_n$ such that $x'\in A^\circ$. Since $x'\leq x_n\leq x$ it follows that $x\in(A^\circ)^\up$. The argument also shows that for each $x\in A$ there exists an $x'\in S^\La_{\rm fin}$ such that $x'\leq x$ and hence $A^\circ\sub S^\La_{\rm fin}$.

We are now ready to prove the statements of the proposition. Clearly $1_{Y^\up}$ is monotone for each $Y\sub S^\La_{\rm fin}$. If $Y$ has finitely many elements, then $1_{Y^\up}$ depends on finitely many coordinates and as a result is clearly continuous. For general $Y\in\Hi(\La)$ we can find $Y_n\in\Hi_{\rm fin}(\La)$ that increase to $Y$. Then $1_{Y^\up_n}\up 1_{Y^\up}$ so $1_{Y^\up}$ is the increasing limit of continuous functions and hence lower semi-continuous. Since $Y\neq\{\un 0\}$ implies $1_{Y^\up}(\un 0)=0$, this shows that the map $Y\mapsto 1_{Y^\up}$ maps $\Hi(\La)$ into $\Li_+(S^\La,\{0,1\})$.

Conversely, if $\phi\in\Li_+(S^\La,\{0,1\})$, then $O_\phi$ is open and increasing so by (\ref{cirup}) $O^\circ_\phi\sub S^\La_{\rm fin}$. Setting $Y:=O^\circ_\phi$ we have $Y^\circ=Y$ by (\ref{elem})~(i) and $Y\neq\{\un 0\}$ by the fact that $\phi(\un 0)=0$. This shows that the map $\phi\mapsto O^\circ_\phi$ maps $\Li_+(S^\La,\{0,1\})$ into $\Hi(\La)$. Using (\ref{elem})~(iii) and (\ref{cirup}) we see that the maps $Y\mapsto 1_{Y^\up}$ and $\phi\mapsto O^\circ_\phi$ are each other's inverses.

Using the fact that a function $f:S^\La\to\{0,1\}$ is continuous with respect to the product topology if and only if it depends on finitely many coordinates \cite[Lemma~4.13]{Swa22}, we see that for $Y\in\Hi(\La)$ one has $1_{Y^\up}\in\Ci_+(S^\La,\{0,1\})$ if and only if $Y\in\Hi_{\rm fin}(\La)$.
\epro

We next start to prepare for the proof of Proposition~\ref{P:topol}. For any metric space $(\Xc,d)$, we let $\Ki_+(\Xc)$ denote the space of all nonempty compact subsets of $\Xc$. The \emph{Hausdorff metric} $d_{\rm H}$ on $\Ki_+(\Xc)$ is defined as
\be\label{Haus}
d_{\rm H}(K_1,K_2):=\sup_{x_1\in K_1}d(x_1,K_2)\vee\sup_{x_2\in K_2}d(x_2,K_1),
\ee
where $d(x,A):=\inf_{y\in A}d(x,y)$ denotes the distance between a point $x\in\Xc$ and a set $A\sub\Xc$. It follows from \cite[Lemma~B.1]{SSS14} that the topology generated by $d_{\rm H}$ only depends on the topology on $\Xc$ and not on the choice of the metric $d$. We call this the \emph{Hausdorff topology}.\footnote{Note the subtle difference between ``the Hausdorff topology'' (the topology generated by the Hausdorff metric) and ``a Hausdorff topology'' (any topology satisfying Hausdorff's separation axiom).} By \cite[Lemma~B.3]{SSS14}, if $\Xc$ is compact, then so is $\Ki_+(\Xc)$.

In analogy with the definition of $A^\up$ in (\ref{circ}), for any $A\sub S^\La$ we set
\be
A^\down:=\big\{x\in S^\La:\exists y\in A\mbox{ s.t.\ }y\geq x\big\}.
\ee
Then a set $A$ is decreasing precisely if $A^\down=A$. We let $A^{\rm c}:=S^\La\beh A$ denote the complement of a set $A\sub S^\La$. We let $\Ki_+(S^\La)$ denote the space of nonempty compact subsets of $S^\La$, equipped with the Hausdorff topology, and set
\be
\Ki^\down_+(S^\La):=\big\{A\in\Ki_+(S^\La):A^\down=A\big\}.
\ee
We equip $\Ki^\down_+(S^\La)$ with the induced topology from $\Ki_+(S^\La)$. The proof of Proposition~\ref{P:topol} is based on the following three lemmas.

\bl[Closed subspace]
The\label{L:dwcls} set $\Ki^\down_+(S^\La)$ is a closed subset of $\Ki_+(S^\La)$.
\el

\bl[Bijection to compact decreasing sets]
The\label{L:bij} map $Y\mapsto(Y^\up)^{\rm c}$ is a bijection from $\Hi(\La)$ to $\Ki^\down_+(S^\La)$.
\el

\bl[Convergence criterion]
For\label{L:duco} $A_n,A\in\Ki^\down_+(S^\La)$, the following statements are equivalent:
\be
{\rm(i)}\ 1_{A_n}(x)\asto{n}1_A(x)\quad\forall x\in S^\La_{\rm fin}
\qquad
{\rm(ii)}\ A_n\asto{n}A\mbox{ in }\Ki^\down_+(S^\La).
\ee
\el

We first show how these lemmas imply Proposition~\ref{P:topol} and then prove the lemmas.\med

\bpro[Proof of Proposition~\ref{P:topol}]
Since $S^\La$ is compact, by \cite[Lemma~B.1]{SSS14} so is $\Ki_+(S^\La)$. Then Lemma~\ref{L:dwcls} implies that also $\Ki^\down_+(S^\La)$ is compact. We can then use the bijection from Lemma~\ref{L:bij} to equip $\Hi(\La)$ with a topology such that $Y_n\to Y$ in $\Hi(\La)$ if and only if $(Y^\up_n)^{\rm c}\to(Y^\up)^{\rm c}$ in $\Ki^\down_+(S^\La)$. Since $\Ki^\down_+(S^\La)$ is compact and metrisable, the same is then true for $\Hi(\La)$. Finally, Lemma~\ref{L:duco} tells us that convergence in $\Hi(\La)$ is equivalent to (\ref{topol}).
\epro

The proof of Lemma~\ref{L:dwcls} needs a bit of preparation. It will be convenient to choose an explicit metric generating the topology on $S^\La$. We choose a bijection $\ga:\La\to\N$ and set
\be\label{metric}
d(x,y):=\sum_{i\in\La}3^{-\ga(i)}1_{\{x(i)\neq y(i)\}}\qquad(x,y\in S^\La).
\ee
Setting $\La_n:=\{i\in\La:\ga(i)\leq n\}$, we observe that
\be\label{3trick}
d(x,y)<3^{-n}\quad\volgt\quad x(i)=y(i)\ \forall i\in\La_n\quad\volgt\quad d(x,y)\leq\ha 3^{-n}.
\ee
Indeed, if $x(i)=y(i)$ for all $i\in\La_n$, then $d(x,y)\leq\sum_{k=n+1}^\infty 3^{-k}=\ha 3^{-n}$, while if $x(i)\neq y(i)$ for some $i\in\La_n$, then $d(x,y)\geq 3^{-n}$.

\bl[Convergence of upset and downset]
Assume\label{L:updownCONV} that $A_n,A\in\Ki_+(S^\La)$ satisfy $A_n\to A$. Then also $A_n^\up\to A^\up$ and $A_n^\down\to A^\down$.
\el

\bpro
By symmetry, it suffices to prove the statement for the upsets. We will show that if $A,B\in\Ki_+(S^\La)$ satisfy $d_{\rm H}(A,B)<3^{-n}$ for some $n\in\N$, then $d_{\rm H}(A^\up,B^\up)<3^{-n}$. It suffices to show that $d(x,B^\up)<3^{-n}$ for each $x\in A^\up$ and $d(y,A^\up)<3^{-n}$ for each $y\in B^\up$. By symmetry, it suffices to prove the latter claim. If $y\in B^\up$, then we can find $b\in B$ such that $b\leq y$. Since $d_{\rm H}(A,B)<3^{-n}$, by (\ref{3trick}) there exists an $a\in A$ such that $a(i)=b(i)$ for all $i\in\La_n$. Defining $x\in S^\La$ by $x(i):=y(i)$ if $i\in\La_n$ and $x(i):=a(i)$ if $i\not\in\La_n$, we see that $a\leq x$ so $x\in A^\up$. By (\ref{3trick}) moreover $d(x,y)<3^{-n}$ so we conclude that $d(y,A^\up)<3^{-n}$.
\epro

\bpro[Proof of Lemma~\ref{L:dwcls}]
Assume that $A_n\in\Ki^\down_+(S^\La)$ satisfy $A_n\to A$ for some $A\in\Ki_+(S^\La)$. We need to show that $A$ is decreasing. Lemma~\ref{L:updownclosed} in the appendix tells us that $A^\down\in\Ki_+(S^\La)$ and Lemma~\ref{L:updownCONV} tells us that $A_n=A_n^\down\to A^\down$. Since $A_n\to A$ and $A_n\to A^\down$ we conclude that $A=A^\down$ so $A$ is decreasing.
\epro

\bpro[Proof of Lemma~\ref{L:bij}]
Let $\Oi^\up_-(S^\La)$ be the set of open increasing sets $A\sub S^\La$ such that $A\neq S^\La$. Then $A\in\Oi^\up_-(S^\La)$ if and only if $1_A\in\Li_+(S^\La,\{0,1\})$ so Proposition~\ref{P:openinc} tells us that $Y\mapsto Y^\up$ is a bijection from $\Hi(\La)$ to $\Oi^\up_-(S^\La)$. The statement now follows from the observation that $A\in\Oi^\up_-(S^\La)$ if and only if $A^{\rm c}\in\Ki^\down_+(S^\La)$.
%Indeed, if $A$ is increasing, $x\not\in A$, and $y\leq x$, then clearly $y\not\in A$. By the same argument the complement of a decreasing set is increasing.
\epro

\bpro[Proof of Lemma~\ref{L:duco}]
It will be convenient to work with the explicit metric in (\ref{metric}). With $\La_n$ as in (\ref{3trick}), we set $S^\La_n:=\{x\in S^\La:x(i)=0\ \forall i\in\La\beh\La_n\}$. We will prove that for $A,B\in\Ki^\down_+(S^\La)$
\be\label{dH3}
d_{\rm H}(A,B)<3^{-n}\quad\desd\quad 1_A(x)=1_B(x)\quad\forall x\in S^\La_n.
\ee
By symmetry, it suffices to prove that
\be
\sup_{a\in A}d(a,B)<3^{-n}\quad\desd\quad 1_A(x)\leq 1_B(x)\quad\forall x\in S^\La_n.
\ee
Assume that $\sup_{a\in A}d(a,B)<3^{-n}$ and $1_A(x)=1$ for some $x\in S^\La_n$. Then there exists an $y\in B$ such that $d(x,y)<3^{-n}$ and hence by (\ref{3trick}) $x(i)=y(i)$ for all $i\in\La_n$, which implies that $x\leq y$ and hence $1_B(x)=1$, proving the implication $\volgt$. To prove the converse implication, assume that $1_A(x)\leq 1_B(x)$ for all $x\in S^\La_n$ and fix $a\in A$. Define $x\in S^\La_n$ by $x(i):=a(i)$ if $i\in\La_n$ and $x(i):=0$ if $i\not\in\La_n$. Then $x\in A$ and hence also $x\in B$. By (\ref{3trick}), this implies that $d(a,B)\leq\ha 3^{-n}$. Since this holds for all $a\in A$, we conclude that $\sup_{a\in A}d(a,B)<3^{-n}$.
\epro

We finally provide the proof of Proposition~\ref{P:ord}.\med

\bpro[Proof of Proposition~\ref{P:ord}]
By (\ref{dualorder}), $Y,Z\in\Hi(\La)$ satisfy $Y\leq Z$ if and only if $Y^\up\sub Z^\up$ which is equivalent to $(Y^\up)^{\rm c}\supset(Z^\up)^{\rm c}$. Recalling the way the topology on $\Hi(\La)$ is defined in the proof of Proposition~\ref{P:topol}, this means that to show that the set in (\ref{ord}) is closed we may equivalently show that
\be
\big\{(A,B)\in\Ki^\down_+(S^\La)^2:A\supset B\big\}
\ee
is a closed subset of $\Ki^\down_+(S^\La)^2$. We need to show that if $A_n\to A$, $B_n\to B$, and $A_n\supset B_n$ for each $n$, then $A\supset B$. This follows from the fact that by \cite[Lemma~B.1]{SSS14}, $A_n\to A$ implies that
\be
A=\big\{x:\exists x_n\in A_n\mbox{ s.t.\ }x_n\to x\big\},
\ee
and similarly for $B$. To complete the proof, we must show that $\emptyset\leq Y\leq Y_{\rm top}$ for all $Y\in\Hi(\La)$. This follows from (\ref{dualorder}) and the observations that $\emptyset^\up=\emptyset$, $Y_{\rm top}^\up=S^\La\beh\{\un 0\}$, and $\un 0\not\in Y$ for all $Y\in\Hi(\La)$ by the definition of $\Hi(\La)$ in (\ref{Hidef}). (Here we have used that if $Y\in\Hi(\La)$ would satisfy $\un 0\in Y$, then by the fact that $Y=Y^\circ$ we would have $Y=\{\un 0\}$ contradicting the definition of $\Hi(\La)$.)
\epro

We conclude this subsection with the following simple observation that will be of use later on.

\bl[Continuity of the duality function]
The\label{L:psico} function $S^\La\ni x\mapsto\psim(x,Y)$ is continuous for each $Y\in\Hi_{\rm fin}(\La)$ and the function $\Hi(\La)\ni Y\mapsto\psim(x,Y)$ is continuous for each $x\in S^\La_{\rm fin}$.
\el

\bpro
If $Y\in\Hi_{\rm fin}(\La)$, then $x\mapsto\psim(x,Y)$ depends on finitely many coordinates and hence is continuous. If $x\in S^\La_{\rm fin}$, then $Y\mapsto\psim(x,Y)$ is continuous by the definition of the topology on $\Hi(\La)$ in (\ref{topol}).
\epro

\subsection{The dual process}\label{S:DualMP}

In this subsection we prove Lemma~\ref{L:lowF}, Theorem~\ref{T:Fell}, and Lemma~\ref{L:Ymon}. The proof of Lemma~\ref{L:lowF} needs two preparatory lemmas.

\bl[Properties of the forward flow]
Assume\label{L:Xmon} (\ref{mas}) and (\ref{downsum}). Then almost surely, for each $s\leq t$ the map $\Xb_{s,t}$ is continuous and monotone with $\Xb_{s,t}(\un 0)=\un 0$. 
\el

\bpro
The maps $\Xb_{s,t}$ are continuous by Theorem~\ref{T:graph}. Using (\ref{mas})~(ii) and Proposition~\ref{P:finap} in the appendix, we see that moreover each map $\Xb_{s,t}$ is monotone with $\Xb_{s,t}(\un 0)=\un 0$. 
\epro

\bl[Backward construction]
Assume\label{L:back} that $m:S^\La\to S^\La$ is continuous and monotone and satisfies $m(\un 0)=\un 0$. Then one has $\phi\circ m\in\Li_+(S^\La,\{0,1\})$ for each $\phi\in\Li_+(S^\La,\{0,1\})$ and $\phi\circ m\in\Ci_+(S^\La,\{0,1\})$ for each $\phi\in\Ci_+(S^\La,\{0,1\})$.
\el

\bpro
For any map $\phi:S^\La\to\{0,1\}$ such that $\phi(\un 0)=0$ one clearly has $\phi\circ m(\un 0)=0$, and if $\phi$ is monotone then so is $\phi\circ m$ by the monotonicity of $m$. If $\phi$ is lower semicontinuous and $x_n\to x$, then $m(x_n)\to m(x)$ by the continuity of $m$ and hence $\limsup_{n\to\infty}\phi\circ m(x_n)\leq\phi\circ m(x)$, proving that $\phi\circ m$ is lower semicontinuous. If $\phi$ is continuous, then clearly so is $\phi\circ m$.
\epro

\bpro[Proof of Lemma~\ref{L:lowF}]
The fact that the map $\Fb_{t,s}$, defined as in (\ref{bdef}), maps the spaces $\Li_+(S^\La,\{0,1\})$ and $\Ci_+(S^\La,\{0,1\})$ into themselves follows from Lemmas \ref{L:Xmon} and \ref{L:back}. The backward stochastic flow property (\ref{bflow}) is immediate from the definitions.
\epro

We will now first prove Lemma~\ref{L:Ymon} and then Theorem~\ref{T:Fell}. Both proofs need some preparations.

\bl[Dual maps]
For\label{L:dualmap} each continuous monotone map $m:S^\La\to S^\La$ that satisfies $m(\un 0)=\un 0$, there exists a unique map $\hat m:\Hi(\La)\to\Hi(\La)$ such that
\be\label{dualmap}
\psim\big(m(x),Y\big)=\psim\big(x,\hat m(Y)\big)\qquad\big(x\in S^\La,\ Y\in\Hi(\La)\big).
\ee
The map $\hat m$ is monotone with respect to the partial order on $\Hi(\La)$ and satisfies $\hat m(\emptyset)=\emptyset$. Moreover, $\hat m$ maps the space $\Hi_{\rm fin}(\La)$ into itself.
\el

\bpro
Filling in the definition of $\psim$, we see that (\ref{dualmap}) is equivalent to
\be\label{dualmap2}
1_{Y^\up}\circ m(x)=1_{\hat m(Y)^\up}(x)\qquad\big(x\in S^\La,\ Y\in\Hi(\La)\big).
\ee
Here $1_{Y^\up}\in\Li_+(S^\La,\{0,1\})$ by Proposition~\ref{P:openinc} and hence $1_{Y^\up}\circ m\in\Li_+(S^\La,\{0,1\})$ by Lemma~\ref{L:back}. Using again Proposition~\ref{P:openinc}, we see that there exists a unique $Z\in\Hi(\La)$ such that $1_{Y^\up}\circ m=1_{Z^\up}$. Setting $\hat m(Y):=Z$ then defines a map $\hat m:\Hi(\La)\to\Hi(\La)$ such that (\ref{dualmap}) holds, and such a map is clearly unique. If $Y\in\Hi_{\rm fin}(\La)$, then Proposition~\ref{P:openinc} tells us that $1_{Y^\up}\in\Ci_+(S^\La,\{0,1\})$, Lemma~\ref{L:back} tells us that $1_{Y^\up}\circ m\in\Ci_+(S^\La,\{0,1\})$, and hence $\hat m(Y)\in\Hi_{\rm fin}(\La)$ by Proposition~\ref{P:openinc}, proving that $\hat m$ maps the space $\Hi_{\rm fin}(\La)$ into itself.

To see that $\hat m$ is monotone with respect to the partial order on $\Hi(\La)$ defined in (\ref{dualorder}), assume that $Y,Z\in\Hi(\La)$ satisfy $Y^\up\sub Z^\up$. Then $1_{Y^\up}\circ m\leq 1_{Z^\up}\circ m$ and hence $1_{\hat m(Y)^\up}\leq 1_{\hat m(Z)^\up}$ proving that $\hat m(Y)^\up\sub\hat m(Z)^\up$. To see that $\hat m(\emptyset)=\emptyset$ it suffices to note that $Y=\emptyset$ implies $1_{Y^\up}=0$ and hence $1_{\hat m(Y)^\up}=0$ which implies $\hat m(Y)=\emptyset$.
\epro

The next lemma says that in order to check (\ref{dualmap}) it suffices to show that it holds for $x\in S^\La_{\rm fin}$.

\bl[Finite configurations]
If\label{L:finsuf} $Y,Z\in\Hi(\La)$ satisfy $\psim(x,Y)=\psim(x,Z)$ for all $x\in S^\La_{\rm fin}$, then $Y=Z$.
\el

\bpro
If $\psim(x,Y)=\psim(x,Z)$ for all $x\in S^\La_{\rm fin}$, then $1_{Y^\up}(x)=1_{Z^\up}(x)$ for all $x\in S^\La_{\rm fin}$ and hence also $1_{(Y^\up)^{\rm c}}(x)=1_{(Z^\up)^{\rm c}}(x)$ for all $x\in S^\La_{\rm fin}$. By (\ref{dH3}) this implies that $(Y^\up)^{\rm c}=(Z^\up)^{\rm c}$. In view of Proposition~\ref{P:openinc}, $Y$ is uniquely determined by $Y^\up$, allowing us to conclude that $Y=Z$.
\epro

Note that Lemma~\ref{L:dualmap} does not say anything about continuity of the dual map $\hat m$. The next lemma fill this gap.

\bl[Continuity of the dual map]
Let\label{L:contdu} $m:S^\La\to S^\La$ be continuous and monotone with $m(\un 0)=\un 0$. Assume that $m$ maps $S^\La_{\rm fin}$ into itself. Then the dual map $\hat m:\Hi(\La)\to\Hi(\La)$ is continuous with respect to the topology on $\Hi(\La)$.
\el

\bpro
We need to show that $Y_n\to Y$ implies $\hat m(Y_n)\to\hat m(Y)$. By (\ref{topol}), this means that we need to show that
\be\label{tbs}
1_{Y_n^\up}(x)\asto{n}1_{Y^\up}(x)\quad\forall x\in S^\La_{\rm fin}
\ee
implies
\be
1_{\hat m(Y_n)^\up}(x)\asto{n}1_{\hat m(Y)^\up}(x)\quad\forall x\in S^\La_{\rm fin}.
\ee
By (\ref{dualmap2}), the latter is equivalent to
\be
1_{Y_n^\up}\circ m(x)\asto{n}1_{Y^\up}\circ m(x)\quad\forall x\in S^\La_{\rm fin}.
\ee
If $m$ maps $S^\La_{\rm fin}$ into itself, then this is indeed implied by (\ref{tbs}).
\epro

\bpro[Proof of Lemma~\ref{L:Ymon}]
Lemma~\ref{L:Xmon} says that for each $s\leq t$ the map $\Xb_{s,t}$ is continuous and monotone with $\Xb_{s,t}(\un 0)=\un 0$. Therefore, by the duality relation (\ref{dual}), in the notation of Lemma~\ref{L:dualmap}, we have
\be\label{YhatX}
\Yb_{t,s}=\hat\Xb_{s,t}\qquad(s\leq t).
\ee
In view of this, the monotonicity of $\Yb_{t,s}$ follows from Lemma~\ref{L:dualmap}.
\epro

\bp[Backward evolution equation]
Under\label{P:backevol} the assumptions of Theorem~\ref{T:Fell}, almost surely, for each $u\in\R$ and $Y\in\Hi(\La)$, there exists a unique cadlag function $(-\infty,u]\ni t\mapsto Y_t\in\Hi(\La)$ such that
\be\label{backevol}
Y_u=Y\quad\mbox{and}\quad Y_{t-}=\left\{\ba{ll}
\dis\hat m(Y_t)\quad&\mbox{if }(m,t)\in\om,\\[5pt]
\dis Y_t\quad&\mbox{otherwise}
\ea\right.\qquad(t\leq u).
\ee
This function is given by $Y_t=\Yb_{u,t}(Y)$ $(t\leq u)$, where $(\Yb_{t,s})_{t\geq s}$ is the backward stochastic flow defined in (\ref{Ybdef}).
\ep

\bpro
Fix $u\in\R$ and $Y\in\Hi(\La)$ and define $(Y_t)_{t\leq u}$ by $Y_t:=\Yb_{u,t}(Y)$ $(t\leq u)$. Then $Y_u=Y$ since $(\Yb_{t,s})_{t\geq s}$ is a backward stochastic flow in the sense of (\ref{bflow}). Fix $x\in S^\La_{\rm fin}$ and set
\be
T:=\big\{t\leq u:\exists(m,t)\in\om\mbox{ s.t.\ }m(x)\neq x\big\}.
\ee
By (\ref{mas}) and (\ref{downsum})~(i), the set $T$ is locally finite, so we can write $T=\{t_k:k\geq 1\}$ with $u\geq t_1\geq t_2\geq\cdots$ and $t_k\to-\infty$ as $k\to\infty$. We let $m_k$ denote the corresponding maps such that $(m_k,t_k)\in\om$. We observe that
\be
\psim(x,Y_s)=\psim\big(x,\Yb_{t,s}(Y_t)\big)=\psim\big(\Xb_{s,t}(x),Y_t\big)\qquad(s\leq t\leq u).
\ee
This shows that $\psim(x,Y_s)=\psim(x,Y_t)$ if $\Xb_{s,t}(x)=x$, so the function $(-\infty,u]\ni t\mapsto\psim(x,Y_t)$ is constant on intervals of the form $(t_{k+1},t_k]$. Since $\Xb_{t_k-\eps,t_k}(x)=m_k(x)$ for $\eps>0$ small enough, we have
\be
\psim(x,Y_{t_k-})=\psi\big(m_k(x),Y_{t_k}\big)=\psi\big(x,\hat m_k(Y_{t_k})\big)\qquad(k\geq 1).
\ee
Since this holds for arbitrary $x\in S^\La_{\rm fin}$, by Lemma~\ref{L:finsuf} and the definition of the topology on $\Hi(\La)$ in (\ref{topol}), it follows that $(Y_t)_{t\leq u}$ is cadlag and satisfies (\ref{backevol}).

To show uniqueness, assume that $(Y_t)_{t\leq u}$ is cadlag and satisfies (\ref{backevol}). Fix $s<u$ and $x\in S^\La_{\rm fin}$ and define $(X_t)_{t\geq s}$ by
\be
X_t:=\Xb_{s,t}(x)\quad(t\geq s).
\ee
We claim that the function
\be\label{pathdu}
[s,u]\ni t\mapsto\psim(X_t,Y_t)
\ee
is constant.

We start by showing that it is cadlag. By Lemma~\ref{L:finsys}, the function $(X_t)_{t\geq s}$ takes values in $S^\La_{\rm fin}$. Combining this with Theorem~\ref{T:graph} we see that it must be piecewise constant and right continuous. Let $s<s_1<s_2<\cdots$ be the times when it jumps, set $s_0:=s$, and $x_k:=X_{s_k}$ $(k\geq 0)$. Then on each interval of the form $[s_{k-1},s_k]$, the function $t\mapsto\psi(x_k,Y_t)$ must be cadlag by the assumption that $(Y_t)_{t\leq u}$ is cadlag and the definition of the topology on $\Hi(\La)$ in (\ref{topol}). This implies that the function in (\ref{pathdu}) is cadlag. Moreover, its left and right limits at a time $t$ are given by $\psim(X_{t-},Y_{t-})$ and $\psim(X_t,Y_t)$, respectively.

We can now use the assumption that $(Y_t)_{t\leq u}$ solves (\ref{backevol}) while $(X_t)_{t\geq s}$ solves (\ref{evo}) to check that
\be
\psim\big(X_{t-},Y_{t-}\big)=\psim\big(X_{t-},\hat m(Y_t)\big)=\psim\big(m(X_{t-}),Y_t\big)=\psim(X_t,Y_t)
\ee
at times $t$ when $(m,t)\in\om$, while more trivialy $\psim(X_{t-},Y_{t-})=\psim(X_t,Y_t)$ at all other times. This shows that the function in (\ref{pathdu}) is continuous. Since it takes values in the set $\{0,1\}$, this implies that it must in fact be constant.

The fact that the function in (\ref{pathdu}) is constant implies that
\be
\psim(x,Y_s)=\psim\big(\Xb_{s,u}(x),Y\big)=\psim\big(x,\Yb_{t,s}(Y)\big).
\ee
Since this holds for all $x\in S^\La_{\rm fin}$ and $s<u$, it follows that solutions to (\ref{backevol}) are unique and given by $Y_t=\Yb_{u,t}(Y)$ $(t\leq u)$.
\epro

\bpro[Proof of Theorem~\ref{T:Fell}]
Using the fact that the backward stochastic flow $(\Fb_{t,s})_{t\geq s}$ is stationary with independent increments, it is straightforward to check that (\ref{Feldef}) defines a Markov process $(Y_t)_{t\geq 0}$ with semigroup $(Q_t)_{t\geq 0}$ defined in (\ref{Qt}). The fact that $(Y_t)_{t\geq 0}$ has caglad sample paths follows from Proposition~\ref{P:backevol}. (Note that because of time reversal, $(Y_t)_{t\geq 0}$ is caglad while the function in (\ref{backevol}) is cadlag.) It therefore remains to prove that $(Q_t)_{t\geq 0}$ is a Feller semigroup. By a well-known characterisation of Feller semigroups \cite[Section~4.2]{Swa22}, letting $\Mi_1(\Hi(\La))$ denote the space of probability measures on $\Hi(\La)$ equipped with the topology of weak convergence, this means that we must show that the map
\be
\Hi(\La)\times\half\ni(Y,t)\mapsto Q_t(Y,\,\cdot\,)\in\Mi_1\big(\Hi(\La)\big)
\ee
is continuous. Since almost sure convergence implies weak convergence in law it suffices to show that for deterministic $(Y_n,t_n)$ and $(Y,t)$,
\be
(Y_n,t_n)\asto{n}(Y,t)\quad\mbox{implies}\quad\Yb_{t_n,0}(Y_n)\asto{n}\Yb_{t,0}(Y)\quad{\rm a.s.}
\ee
By the definition of the topology on $\Hi(\La)$ in (\ref{topol}), we must show that
\be
1_{\Yb_{t_n,0}(Y_n)^\up}(x)\asto{n}1_{\Yb_{t,0}(Y)^\up}(x)\quad{\rm a.s.}\quad(x\in S^\La_{\rm fin}).
\ee
By duality (\ref{dual}), this is equivalent to
\be\label{Ytco}
1_{Y_n^\up}\big(\Xb_{0,t_n}(x)\big)\asto{n}1_{Y^\up}\big(\Xb_{0,t}(x)\big)\quad{\rm a.s.}\quad(x\in S^\La_{\rm fin}).
\ee
By Lemma~\ref{L:finsys}, the process $t\mapsto\Xb_{0,t_n}(x)$ is a continuous-time Markov chain with countable state space $S^\La_{\rm fin}$. Since $t$ is deterministic, it is a.s.\ not a jump time of this continuous-time Markov chain, so $\Xb_{0,t_n}(x)=\Xb_{0,t}(x)$ for all $n$ large enough. Using the fact that $\Xb_{0,t}(x)\in S^\La_{\rm fin}$, the convergence in (\ref{Ytco}) then follows from the definition of the topology on $\Hi(\La)$ in (\ref{topol}).
\epro

\subsection{The upper invariant law}\label{S:UP&Surv}

In this subsection we prove Theorem~\ref{T:upp} and Lemma~\ref{L:surstab}.\med

\bpro[Proof of Theorem~\ref{T:upp}]
Immediate from Lemma~\ref{L:Ymon}, Proposition~\ref{P:ord}, and Proposition~\ref{P:upp} in the appendix.
\epro

\bpro[Proof of Lemma~\ref{L:surstab}]
To prove part~(i), we observe using duality (\ref{dual}) that for each $Y\in\Hi_{\rm fin}(\La)$,
\be\ba{l}
\dis\E\big[\psim\big(\Xb_{0,t}(\un\top),Y\big)\big]
=\E\big[\psim\big(\un\top,\Yb_{t,0}(Y)\big)\big]\\[5pt]
\dis\quad=\P\big[\Yb_{0,-t}(Y)\neq\emptyset\big]
\asto{t}\P\big[\Yb_{0,-t}(Y)\neq\emptyset\ \forall t\geq 0\big],
\ec
so taking the limit, using (\ref{Xup}) and Lemma~\ref{L:psico}, we see that
\be\label{nuXY}
\int\ov\nu_{\rm X}(\di x)\psim(x,Y)=\P\big[\Yb_{0,-t}(Y)\neq\emptyset\ \forall t\geq 0\big]\qquad\big(Y\in\Hi_{\rm fin}(\La)\big).
\ee
In particular
\be
\ov\nu_{\rm X}\big(\{x\in S^\La:x\geq y\}\big)=\P\big[\Yb_{0,-t}(\{y\})\neq\emptyset\ \forall t\geq 0\big]\qquad(y\in S^\La_{\rm fin}\beh\{\un 0\}).
\ee
If $(Y_t)_{t\geq 0}$ dies out, then this is zero for all $y\in S^\La_{\rm fin}\beh\{\un 0\}$, proving that $\ov\nu_{\rm X}=\de_{\un 0}$. Conversely, if $\ov\nu_{\rm X}=\de_{\un 0}$, then the left-hand side of (\ref{nuXY}) is zero for all $Y\in\Hi_{\rm fin}(\La)$ so $(Y_t)_{t\geq 0}$ dies out.

The proof of part~(ii) is similar. Using duality (\ref{dual}) we have, for each $x\in S^\La_{\rm fin}$,
\be\ba{l}
\dis\E\big[\psim\big(x,\Yb_{0,-t}(Y_{\rm top})\big)\big]
=\E\big[\psim\big(\Xb_{-t,0}(x),Y_{\rm top}\big)\big]\\[5pt]
\dis\quad=\P\big[\Xb_{0,t}(x)\neq\un 0\big]
\asto{t}\P\big[\Xb_{0,t}(x)\neq\un 0\ \forall t\geq 0\big],
\ec
so taking the limit, using (\ref{Yup}) and Lemma~\ref{L:psico}, we see that
\be\label{nuYX}
\int\ov\nu_{\rm Y}(\di Y)\psim(x,Y)=\P\big[\Xb_{0,t}(x)\neq\un 0\ \forall t\geq 0\big]\qquad(x\in S^\La_{\rm fin}).
\ee
Let $\ov Y$ be a random variable with law $\ov\nu_{\rm Y}$. If $\ov Y=\emptyset$ a.s., then the quantity in (\ref{nuYX}) is clearly zero for all $x\in S^\La_{\rm fin}$. On the other hand, if $\ov Y\neq\emptyset$ with positive probability, then by the fact that $Y$ is a random subset of $S^\La_{\rm fin}$ there must exist an $x\in S^\La_{\rm fin}$ for which the quantity in (\ref{nuYX}) is positive.
\epro

\subsection{Homogeneous invariant laws}\label{S:homlaw}

In this subsection we prepare for the proof of Theorem~\ref{T:ergo} by proving a more general statement that, specialised to the cooperative contact process, will yield Theorem~\ref{T:ergo}. We first need a characterisation of the upper invariant law of the dual process $(Y_t)_{t\geq 0}$ in terms of the forward process $(X_t)_{t\geq 0}$. We work for the moment in the general set-up of Theorem~\ref{T:upp}.

\bl[Distribution determining functions]
Let\label{L:disdet} $Y,Z$ be $\Hi(\La)$-valued random variables such that
\be
\E\big[\prod_{k=1}^n\psim(x_k,Y)\big]=\E\big[\prod_{k=1}^n\psim(x_k,Z)\big]
\qquad\forall n\geq 1,\ x_1,\ldots,x_n\in S^\La_{\rm fin}.
\ee
Then $Y$ and $Z$ are equal in law.
\el

\bpro
Let $\Fi$ be the class of functions $f:\Hi(\La)\to\R$ of the form $f(Y)=\prod_{k=1}^n\psim(x_k,Y)$ with $n\geq 1$ and $x_1,\ldots,x_n\in S^\La_{\rm fin}$. By \cite[Lemma~4.37]{Swa22} it suffices to show that each $f\in\Fi$ is continuous, and that $\Fi$ is closed under products and separates points in the sense that for each $Y,Z\in\Hi(\La)$, there exists an $f\in\Fi$ such that $f(Y)\neq f(Z)$. Functions $f\in\Fi$ are continuous by Lemma~\ref{L:psico}, the class $\Fi$ is closed under products by construction, and $\Fi$ separates points by Lemma~\ref{L:finsuf}.
\epro

\bl[Characterisation of the upper invariant law]
Assume\label{L:upchar} (\ref{mas}), and assume that the rates $(r_m)_{m\in\Gi}$ satisfy (\ref{downsum}) and (\ref{upsum}). Then the upper invariant law $\ov\nu_{\rm Y}$ of the dual process, defined in (\ref{Yup}) is uniquely characterised by the fact that
\be\label{upchar}
\int\ov\nu_{\rm Y}(\di Y)\prod_{k=1}^n\psim(x_k,Y)=\P\big[\Xb_{0,t}(x_k)\neq\un 0\ \forall t\geq 0,\ 1\leq k\leq n\big]
\qquad(x_1,\ldots,x_n\in S^\La_{\rm fin}).
\ee
\el

\bpro
For $n=1$ formula (\ref{upchar}) has already been proved as formula (\ref{nuYX}). The proof for general $n\geq 1$ is completely the same. The fact that $\ov\nu_{\rm Y}$ is uniquely characterised by (\ref{upchar}) follows from Lemma~\ref{L:disdet}.
\epro

We now set out to prove a general result in the spirit of Theorem~\ref{T:ergo}. We assume from now on that the grid $\La$ is a (not necessarily abelian) group with product $(i,j)\mapsto ij$. A bit unusually, to stay closer to the notation for $\Z^d$, we will denote the unit element of $\La$ by $0$. For each $i\in\La$ we define shift operators $S_i$ and $T_i$ acting on configurations $x\in S^\La$ and local maps $m:S^\La\to S^\La$, respectively, by
\be\label{shift}
(S_ix)(j):=x(i^{-1}j)\quad(j\in\La)\quand(T_im)(x)(j):=m(S_{i^{-1}}x)(i^{-1}j)\quad(j\in\La).
\ee
We say that a probability measure $\mu$ on $S^\La$ is \emph{homogeneous} if $\mu=\mu\circ S_i^{-1}$ for all $i\in\La$, and we say that the rates $(r_m)_{m\in\Gi}$ are \emph{translation invariant} if 
\be\label{trainv}
T_im\in\Gi\quand r_{T_im}=r_m\quad(i\in\La,\ m\in\Gi).
\ee
As we will show, Theorem~\ref{T:ergo} is a simple consequence of the following more general theorem. Recall from (\ref{xabs}) that $|x|:=\sum_{i\in\La}1_{\{x(i)\neq 0\}}$. The assumptions $S=\{0,\ldots,n\}$ and $|\Di(m)|=1$ $(m\in\Gi)$ below are probably not needed but they significantly simplify the proof. The theorem applies to the cooperative contact process and more generally to spin systems, which is more than sufficient for our purposes.

\bt[Homogeneous initial laws]
Assume\label{T:ergo2} (\ref{mas}) and that the grid $\La$ is a group. Assume that $S=\{0,\ldots,n\}$ and $|\Di(m)|=1$ for all $m\in\Gi$. Assume that the rates $(r_m)_{m\in\Gi}$ satisfy (\ref{downsum}) and (\ref{upsum}) and are translation invariant. Assume moreover that:
\begin{itemize}
\item[{\rm(i)}] $\dis\forall\eps>0\ \exists N<\infty$ such that $|x|\leq N$ implies $\dis\P^x\big[\exists t\geq 0\mbox{ s.t.\ }X_t=\un 0\big]\geq\eps$.
\item[{\rm(ii)}] $\dis\P^x\big[X_t\geq y\big]>0$ for each $x,y\in S^\La_{\rm fin}\beh\{\un 0\}$ and $t>0$.
\end{itemize}
Then the monotone dual $(Y_t)_{t\geq 0}$ started in an initial law that is homogeneous with $Y_0\neq\emptyset$ a.s.\ satisfies
\be\label{ergo2}
\P\big[Y_t\in\,\cdot\,\big]\Asto{t}\ov\nu_{\rm Y},
\ee
where $\Rightarrow$ denotes weak convergence of probability measures on $\Hi(\La)$.
\et

The basic idea behind the proof of Theorem~\ref{T:ergo2} is very old and goes back to the work of Vasil'ev \cite{Vas69} and Harris \cite{Har76}. The details differ, however, from model to model and depend on the type of duality that is being used. We first prove two lemmas and then prove the theorem.

\bl[Extinction versus unbounded growth]
Assume\label{L:exgro} (\ref{mas}), (\ref{downsum}), (\ref{upsum}) and condition~(i) of Theorem~\ref{T:ergo2}. Then
\be\label{exgro}
\P^x\big[0<|X_t|<N\big]\asto{t}0\qquad(x\in S^\La_{\rm fin},\ N<\infty).
\ee
\el

\bpro
Condition~(i) of Theorem~\ref{T:ergo2} says that each time the process returns to a state $x$ with $|x|\leq N$, there is a probability of at least $\eps$ that the process gets extinct. By a standard argument, this implies that almost surely either $X_t=\un 0$ for some $t\geq 0$ or $|X_t|\to\infty$, see \cite[Lemma~6.36]{Swa22}.
\epro

The next lemma is the key step in the proof of Theorem~\ref{T:ergo2}.

\bl[Large is good]
Assume\label{L:large} (\ref{mas}) and that the grid $\La$ is a group. Assume that the rates $(r_m)_{m\in\Gi}$ satisfy (\ref{downsum}) and (\ref{upsum}) and are translation invariant, and assume condition~(ii) of Theorem~\ref{T:ergo2}. Let $(Y_t)_{t\geq 0}$ be the monotone dual process, started in a homogeneous initial law with $Y_0\neq\emptyset$ a.s. Then for each $s,\eps>0$, there exists an $N<\infty$ such that for any $x\in S^\La_{\rm fin}$
\be
|x|\geq N\quad\mbox{implies}\quad\P\big[\psim(x,Y_s)=0\big]\leq\eps.
\ee
\el

\bpro
We construct $(Y_t)_{t\geq 0}$ as $Y_t:=\Yb_{s,s-t}(Y_0)$ $(t\geq 0)$ where $Y_0$ is independent of the graphical representation $\om$, and use the duality relation (\ref{dual}) to write
\be\label{our}
\P\big[\psim(x,Y_s)=0\big]=\P\big[\psim(\Xb_{0,s}(x),Y_0)=0\big]=\int\P\big[Y_0\in\di Y\big]\P\big[1_{Y^\up}\big(\Xb_{0,s}(x)\big)=0\big].
\ee
Fix $s,\eps>0$, $N<\infty$, and $x\in S^\La$ with $|x|\geq N$. Let $\De_k\sub\La$ be finite sets such that $0\in\De_k$ and $\De_k\up\La$ as $k\to\infty$. For $i\in\La$ write $i\De_k:=\{ij:j\in\De_k\}$. It is easy to see that there exists a finite set $\Ga_k\sub\La$ with $|\Ga_k|\geq N/|\De_k|$ such that
\be
x(i)\neq 0\ \forall i\in\Ga_k\mbox{ and the sets }\{i\De_k:i\in\Ga_k\}\mbox{ are disjoint.}
\ee
For each $i\in\La$, let $\om^{i,k}:=\{(m,t)\in\om:\Di(m)\sub i\De_k\}$ and let $(\Xb^{i,k}_{s,t})_{s\leq t}$ be the stochastic flow defined as in (\ref{evo}) but with $\om$ replaced by $\om^{i,k}$. Set
\be
X_t:=\Xb_{0,t}(x)\quand X^{i,k}_t:=\Xb^{i,k}_{0,t}(e_i)\qquad(t\geq 0,\ i\in\Ga_k).
\ee
Then $X^{i,k}_t(j)=0$ for all $j\not\in i\De_k$ and $t\geq 0$, and using the assumptions that $S=\{0,\ldots,n\}$ and $|\Di(m)|=1$ $(m\in\Gi)$ and Proposition~\ref{P:finap} in the appendix it is easy to see that $X_t\geq X^{i,k}_t$ for each $t\geq 0$. Moreover, since the sets $i\De_k$ are disjoint, the processes $(X^{i,k}_t)_{t\geq 0}$ with $i\in\Ga_k$ are independent. We note that this is the only place in the proof of Theorem~\ref{T:ergo2} where the assumptions $S=\{0,\ldots,n\}$ and $|\Di(m)|=1$ $(m\in\Gi)$ are used. We can estimate the right-hand side of (\ref{our}) as
\be\ba{l}
\dis\int\P\big[Y_0\in\di Y\big]\P\big[1_{Y^\up}(X_s)=0\big]\leq\int\P\big[Y_0\in\di Y\big]\prod_{i\in\Ga_k}\P\big[1_{Y^\up}(X^{i,k}_s)=0\big]\\[5pt]
\dis\quad\stackrel{1}{\leq}\prod_{i\in\Ga_k}\Big(\int\P\big[Y_0\in\di Y\big]\P\big[1_{Y^\up}(X^{i,k}_s)=0\big]^{|\Ga_k|}\Big)^{1/|\Ga_k|}\\[5pt]
\dis\quad\stackrel{2}{=}\prod_{i\in\Ga_k}\Big(\int\P\big[Y_0\in\di Y\big]\P\big[1_{Y^\up}(X^{0,k}_s)=0\big]^{|\Ga_k|}\Big)^{1/|\Ga_k|}\\[5pt]
\dis\quad=\int\P\big[Y_0\in\di Y\big]\P\big[1_{Y^\up}(X^{0,k}_s)=0\big]^{|\Ga_k|}\leq\int\P\big[Y_0\in\di Y\big]\P\big[1_{Y^\up}(X^{0,k}_s)=0\big]^{N/|\De_k|},
\ec
where the inequality $\stackrel{1}{\leq}$ follows from H\"older's inequality and in the equality $\stackrel{2}{=}$ we have used the homogeneity of the law of $Y_0$. Set
\be
Z_k:=\big\{y\in S^\La_{\rm fin}:\P\big[X^{0,k}_s\geq y\big]>0\big\}.
\ee
Since $\De_k\up\La$, it is not hard to see that $X^{0,k}_s$ increases to $\Xb_{0,s}(e_0)$ as $k\to\infty$. As a consequence,
\be
Z_k\up Z:=\big\{y\in S^\La_{\rm fin}:\P\big[\Xb_{0,s}(e_0)\geq y\big]>0\big\}=S^\La_{\rm fin}
\ee
where the final equality follows from condition~(ii) of Theorem~\ref{T:ergo2}. Using this and the assumption that $Y_0\neq\emptyset$ a.s., we see that
\be
\int\P\big[Y_0\in\di Y\big]\P\big[1_{Y^\up}(X^{0,k}_s)=0\big]^{N/|\De_k|}\asto{N}\P\big[Y_0\cap Z_k=\emptyset\big]\asto{k}\P\big[Y_0=\emptyset\big]=0.
\ee
In particular, for each $\eps>0$ we can first choose $k$ large enough and then $N$ large enough so that this expression is $\leq\eps$, which by our previous calculations proves the claim of the lemma.
\epro

\bpro[Proof of Theorem~\ref{T:ergo2}]
Since the space $\Hi(\La)$ is compact, the laws $(\P[Y_t\in\,\cdot\,])_{t\geq 0}$ are tight, so it suffices to prove that $\ov\nu_{\rm Y}$ is the only cluster point. By Lemma~\ref{L:upchar} (using also Lemma~\ref{L:psico}), it suffices to show that for each $n\geq 1$ and $x_1,\ldots,x_n\in S^\La_{\rm fin}$,
\be\label{nuaim}
\E\big[\prod_{k=1}^n\psim(x_k,Y_t)\big]\asto{t}\P\big[\Xb_{0,t}(x_k)\neq\un 0\ \forall t\geq 0,\ 1\leq k\leq n\big]=:\rho(x_1,\ldots,x_n).
\ee
We construct $(Y_t)_{t\geq 0}$ as $Y_t=\Yb_{t,0}(Y_0)$ $(t\geq 0)$ where $Y_0$ is independent of the graphical representation, fix $s>0$, and use duality to rewrite the left-hand side of (\ref{nuaim}) as
\be
\E\big[\prod_{k=1}^n\psim\big(\Xb_{0,t-s}(x_k),\Yb_{t,t-s}(Y_0)\big)\big].
\ee
Here $\Yb_{t,t-s}(Y_0)$ is equally distributed with $Y_s$. Since disjoint parts of the graphical representation are independent, it is independent of $(\Xb_{0,t-s}(x_k))_{1\leq k\leq n}$. If $t$ tends to infinity then so does $t-s$, so we see that to prove (\ref{nuaim}), it suffices to prove that
\be\label{nuaim2}
\E\big[\prod_{k=1}^n\psim\big(\Xb_{0,t}(x_k),Y_s\big)\big]\asto{t}\rho(x_1,\ldots,x_n),
\ee
where $Y_s$ is independent of $(\Xb_{0,t}(x_k))_{1\leq k\leq n}$. We fix $\eps>0$, we choose $N$ in dependence on $s$ and $\eps$ as in Lemma~\ref{L:large}, and introduce the events
\bc
\dis A&:=&\dis\big\{\Xb_{0,t}(x_k)=\un 0\mbox{ for some }1\leq k\leq n\big\},\\[5pt]
\dis B&:=&\dis\big\{0<|\Xb_{0,t}(x_k)|<N\mbox{ for some }1\leq k\leq n\big\},\\[5pt]
\dis C&:=&\dis\big\{|\Xb_{0,t}(x_k)|\geq N\mbox{ for all }1\leq k\leq n\big\}.
\ec
We rewrite the left-hand side of (\ref{nuaim2}) as
\be\ba{l}
\dis\E\big[\prod_{k=1}^n\psim\big(\Xb_{0,t}(x_k),Y_s\big)\,\big|\,A\big]\cdot\P[A]+\E\big[\prod_{k=1}^n\psim\big(\Xb_{0,t}(x_k),Y_s\big)\,\big|\,B\big]\cdot\P[B]\\[5pt]
\dis+\E\big[\prod_{k=1}^n\psim\big(\Xb_{0,t}(x_k),Y_s\big)\,\big|\,C\big]\cdot\P[C].
\ec
Here the first term is zero since $\psim(x,Y)=0$ if $x=\un 0$ and the second tends to zero by Lemma~\ref{L:exgro}. By Lemma~\ref{L:large} and a simple union bound,
\be
\P\big[\psim\big(\Xb_{0,t}(x_k),Y_s\big)=0\mbox{ for some }1\leq k\leq n\,\big|\,C\big]\leq n\eps,
\ee
while by Lemma~\ref{L:exgro} the probability of the event $C$ tends to $\rho(x_1,\ldots,x_n)$. It follows that
\be\ba{r@{\,}l}
\dis(1-n\eps)\rho(x_1,\ldots,x_n)&\dis\leq\liminf_{t\to\infty}\E\big[\prod_{k=1}^n\psim\big(\Xb_{0,t}(x_k),Y_s\big)\big]\\[5pt]
&\dis\leq\limsup_{t\to\infty}\E\big[\prod_{k=1}^n\psim\big(\Xb_{0,t}(x_k),Y_s\big)\big]\leq\rho(x_1,\ldots,x_n).
\ec
Since $\eps>0$ is arbitrary, this proves the theorem.
\epro

\subsection{The cooperative contact process}\label{S:coop}

In this subsection we prove Theorems \ref{T:ergo} and \ref{T:tetcont}. We will show that most statements actually remain true if the grid $\Z^d$ is replaced by a general Cayley graph. Throughout this subsection we assume that $\La$ is a finitely generated group and that $\De$ is a finite subset of $\La$ that does not contain $0$, is symmetric in the sense that $k\in\De$ implies $k^{-1}\in\De$, and that generates $\La$. To avoid trivialities we also assume $|\De|\geq 2$. We equip $\La$ with the structure of a locally finite graph with set of edges
\be\label{Edges}
E:=\big\{\{j,jk\}:j\in\La,\ k\in\De\big\}.
\ee
This says that $(\La,E)$ is the Cayley graph associated with $\La$ and $\De$. Note that if $\{j,k\}\in E$ and $i\in\La$, then $\{ij,ik\}\in E$, and as a result the cooperative contact process on the graph $(\La,E)$ has translation invariant rates in the sense of (\ref{trainv}). Setting $\La:=\Z^d$ with the usual additive group structure and $\De:=\{k:|k|=1\}$ yields the nearest-neighbour process on $\Z^d$.\med

\bpro[Proof of Theorem~\ref{T:ergo}]
Most of the statement remains true if $\Z^d$ with nearest neighbour edges is replaced by a general Cayley graph as described above. We first prove the statement under the assumptions $\de>0$ and $\al<1$. These assumptions imply that conditions (i) and (ii), respectively, of Theorem~\ref{T:ergo2} are satisfied, so in this case formula (\ref{ergo}) follows from (\ref{ergo2}).

We next consider the case $\de=0$ and $\al<1$. Since $\de=0$, we have $\P^x[X_t\neq\un 0\ \forall t\geq 0]=1$ for all $x\neq\un 0$, so by (\ref{nuaim}) to prove the claim it suffices to show that
\be
\E\big[\prod_{k=1}^n\psim(x_k,Y_t)\big]\asto{t}1
\ee
for all $x_1,\ldots,x_n\in\{0,1\}^\La_{\rm fin}\beh\{\un 0\}$. It suffices to show this for $n=1$, so using the duality relation (\ref{dual}), we need to show that
\be
\E\big[\psim\big(\Xb_{0,t}(x),Y_0\big)\big]\asto{t}1\qquad\big(x\in\{0,1\}^\La_{\rm fin}\beh\{\un 0\}\big).
\ee
Filling in the definition of $\psim$, this says that
\be\label{denul}
\P\big[\exists y\in Y_0\mbox{ s.t.\ }\Xb_{0,t}(x)\geq y\big]\asto{t}1\qquad\big(x\in\{0,1\}^\La_{\rm fin}\beh\{\un 0\}\big).
\ee
Since $\de=0$ and $\al<1$ and since $\De$ generates $\La$, it is easy to see that as a consequence of branching, $\Xb_{0,t}(x)\up\un 1$ a.s.\ as $t\to\infty$. Therefore, using the assumption $Y_0\neq\emptyset$ a.s., we see that (\ref{denul}) holds.

We finally consider the case $\de>0$ and $\al=1$. Only in this case we use the assumption that the grid is $\Z^d$ with nearest-neighbour edges. Together with $\al=1$, this has the effect that the cooperative contact process, started in a finite interval, rectangle, or (hyper)cube cannot escape from such a set and hence a.s.\ dies out. This means that the right-hand side of (\ref{nuaim}) is zero for all $x_1,\ldots,x_n\in\{0,1\}^{\Z^d}_{\rm fin}$ and (\ref{nuaim}) is trivially satisfied.
\epro

We next start to prepare for the proof of Theorem~\ref{T:tetcont}. We consider cooperative contact processes on general Cayley graphs as explained at the beginning of this subsection.

\bl[Limits of invariant laws]
Let\label{L:liminv} $\Ii_{\rm Y}(\al,\de)$ denote the set of invariant laws of the monotone dual of the cooperative contact process with parameters $\al$ and $\de$. Assume that $(\al_n,\de_n)\in[0,1]\times\half$ converge to a limit $(\al,\de)$ and that $\nu_n\in\Ii_{\rm Y}(\al_n,\de_n)$ converge weakly on $\Hi(\La)$ to a probility law $\nu$ on $\Hi(\La)$. Then $\nu\in\Ii_{\rm Y}(\al,\de)$.
\el

\bpro
Let $\om^n,\om$ be graphical representations corresponding to the rates $(\al_n,\de_n)$ and $(\al,\de)$, respectively. For each $n$, let $Y^n_0$ have law $\nu_n$ and be independent of $\om^n$. Likewise, let $Y_0$ have law $\nu$ and be independent of $\om$. To show that $\nu$ is invariant, we will show that $\Yb_{t,0}(Y_0)$ has law $\nu$ for each $t\geq 0$. Fix $t\geq 0$. We will show that $\Yb^n_{t,0}(Y^n_0)$ converges weakly in law to $\Yb_{t,0}(Y_0)$. Since each $\nu_n$ is invariant and $\nu_n\Rightarrow\nu$ this then implies that $\Yb_{t,0}(Y_0)$ has law $\nu$. It suffices to show that for a suitable coupling $\Yb^n_{t,0}(Y^n_0)$ converges a.s.\ to $\Yb_{t,0}(Y_0)$. By the definition of the topology on $\Hi(\La)$ in (\ref{topol}) this amounts to showing that
\be
\psim\big(x,\Yb^n_{t,0}(Y^n_0)\big)\asto{n}\psim\big(x,\Yb_{t,0}(Y_0)\big)\quad{\rm a.s.}\qquad(x\in S^\La_{\rm fin}).
\ee
By the duality relation (\ref{dual}) this is equivalent to
\be\label{coupco}
\psim\big(\Xb^n_{0,t}(x),Y^n_0\big)\asto{n}\psim\big(\Xb_{0,t}(x),Y_0\big)\quad{\rm a.s.}\qquad(x\in S^\La_{\rm fin}).
\ee
By Lemma~\ref{L:rateconv} in the appendix, we can couple the graphical representations $\om^n,\om$ in such a way that for each $x\in S^\La_{\rm fin}$ there exists an $N<\infty$ such that $\Xb^n_{0,t}(x)=\Xb_{0,t}(x)$ for all $n\geq N$. By Skorohod's representation theorem \cite[Cor~3.1.6 and Thm~3.1.8]{EK86}, we can couple the $Y^n_0,Y_0$, which are independent of everything else, in such a way that $Y^n_0\to Y_0$ a.s. Then (\ref{coupco}) follows from Lemma~\ref{L:psico}.
\epro

We continue to consider cooperative contact processes on general Cayley graphs. We let $\ov\nu_{\rm X}(\al,\de)$ and $\ov\nu_{\rm Y}(\al,\de)$ denote the upper invariant laws of the process and its monotone dual, in dependence on the parameters $\al$ and $\de$.

\bl[Monotone coupling]
Assume\label{L:monnu} that $\al\leq\al'$ and $\de\leq\de'$. Then $\tet(\al,\de)\geq\tet(\al',\de')$, $\ov\nu_{\rm X}(\al,\de)\geq\ov\nu_{\rm X}(\al',\de')$ in the stochastic order on $S^\La$, and $\ov\nu_{\rm Y}(\al,\de)\geq\ov\nu_{\rm Y}(\al',\de')$ in the stochastic order on $\Hi(\La)$.
\el

\bpro
Let $R_1:=|\Ni_i|$ and $R_2:=|\Ni^{(2)}_i|$ which by translation invariance do not depend on $i\in\La$. Let
\be\ba{c}
\dis\Gi_{\rm dth}:=\big\{{\tt dth}_j:j\in\La\big\},\quad
\Gi_{\rm bra}:=\big\{{\tt bra}_{ij}:j\in\La,\ i\in\Ni_j\big\},\\[5pt]
\dis\Gi_{\rm coop}:=\big\{{\tt coop}_{ii'j}:j\in\La,\ (i,i')\in\Ni^{(2)}_j\big\}.
\ec
Let $\om_{\rm dth}$ and $\om'_{\rm dth}$ be Poisson point sets on $\Gi_{\rm dth}\times\R$ with intensities $\de$ and $\de'-\de$, let $\om_{\rm bra}$ be a Poisson point set on $\Gi_{\rm bra}\times\R$ with intensity $(1-\al')/R_1$, and let $\om_{\rm coop}$ and $\om'_{\rm coop}$ be Poisson point sets on $\Gi_{\rm coop}\times\R$ with intensities $\al/R_2$ and $(\al'-\al)/R_2$. Assume that all these Poisson point sets are independent. Define
\be
\om'_{\rm bra}:=\big\{({\tt bra}_{ij},t):({\tt coop}_{ii'j},t)\in\om'_{\rm coop}\big\}.
\ee
Then setting
\be
\om:=\om_{\rm dth}\cup\om_{\rm bra}\cup\om'_{\rm bra}\cup\om_{\rm coop}
\quand
\om':=\om_{\rm dth}\cup\om'_{\rm dth}\cup\om_{\rm bra}\cup\om_{\rm coop}\cup\om'_{\rm coop}
\ee
defines two graphical representations for two cooperative contact processes, the first one with parameters $\al$ and $\de$ and the second one with parameters $\al'$ and $\de'$. Note that in the latter compared to the former, the death map is applied more frequently while some applications of a branching map ${\tt bra}_{ij}$ have been replaced by an application of a cooperative branching map ${\tt coop}_{ii'j}$ (with the same $i$ and $j$).

Let $(\Xb_{s,t})_{s\leq t}$ and $(\Xb'_{s,t})_{s\leq t}$ be the stochastic flows constructed from $\om$ and $\om'$. Using Proposition~\ref{P:finap} in the appendix, it is straightforward to check that
\be\label{foror}
x\geq x'\quad\mbox{implies}\quad\Xb_{s,t}(x)\geq\Xb'_{s,t}(x')\qquad\big(s\leq t,\ x,x'\in\{0,1\}^\La\big).
\ee
Applying this with $x=x'=e_0$ shows that $\tet(\al,\de)\geq\tet(\al',\de')$ while setting $x=x':=\un 1$, using (\ref{Xup}) we see that $\ov\nu_{\rm X}(\al,\de)\geq\ov\nu_{\rm X}(\al',\de')$ in the stochastic order on $S^\La$. To show that $\ov\nu_{\rm Y}(\al,\de)\geq\ov\nu_{\rm Y}(\al',\de')$ in the stochastic order on $\Hi(\La)$, by (\ref{Yup}) it suffices to show that
\be
\Yb_{t,s}(Y_{\rm top})\geq\Yb'_{t,s}(Y_{\rm top})\quad(t\geq s)
\ee
in the stochastic order on $\Hi(\La)$. By (\ref{dualorder}) this is equivalent to
\be
1_{\Yb_{t,s}(Y_{\rm top})^\up}(x)\geq 1_{\Yb_{t,s}(Y'_{\rm top})^\up}(x)\quad(t\geq s,\ x\in\{0,1\}^\La).
\ee
Here, by the duality relation (\ref{dual}),
\be
1_{\Yb_{t,s}(Y_{\rm top})^\up}(x)=\psim\big(x,\Yb_{t,s}(Y_{\rm top})\big)=\psim\big(\Xb_{s,t}(x),Y_{\rm top}\big)=1_{\{\Xb_{s,t}(x)\neq\un 0\}}
\ee
and likewise for the stochastic flow defined by $\om'$, so we need to show, for each $t\geq s$ and $x\in\{0,1\}^\La$, that $\Xb'_{s,t}(x)\neq\un 0$ implies $\Xb_{s,t}(x)\neq\un 0$. This follows from (\ref{foror}), so the proof is complete.
\epro

\bl[Increasing limit]
Assume\label{L:inclim} that $\al_n,\al\in[0,1]$ and $\de_n,\de\in\half$ satisfy $\al_n\up\al$ and $\de_n\up\de$. Then
\be
\ov\nu_{\rm Y}(\al_n,\de_n)\Asto{n}\ov\nu_{\rm Y}(\al,\de)
\ee
where $\Rightarrow$ denotes weak convergence of probability measures on $\Hi(\La)$.
\el

\bpro
Since $\Hi(\La)$ is compact the measures $\ov\nu_{\rm Y}(\al_n,\de_n)$ are tight so it suffices to show that $\ov\nu_{\rm Y}(\al,\de)$ is their only cluster point. By Lemma~\ref{L:liminv} each cluster point $\nu$ is an invariant law and hence $\nu\leq\ov\nu_{\rm Y}(\al,\de)$. If this is not an equality then there exists a continuous monotone function $f$ such that
\be
\int f\,\di\ov\nu_{\rm Y}(\al_n,\de_n)<\int f\,\di\ov\nu_{\rm Y}(\al,\de)
\ee
for $n$ large enough, contradicting Lemma~\ref{L:monnu}, so we conclude that $\nu=\ov\nu_{\rm Y}(\al,\de)$.
\epro

\bl[Decreasing limit]
Let\label{L:declim} $A$ be the set defined in (\ref{alde}). Assume that $(\al_n,\de_n)\in A$ satisfy $\al_n\down\al$ and $\de_n\down\de$ for some $(\al,\de)\in A$. Then
\be
\ov\nu_{\rm Y}(\al_n,\de_n)\Asto{n}\ov\nu_{\rm Y}(\al,\de).
\ee
%where $\Rightarrow$ denotes weak convergence of probability measures on $\Hi(\La)$.
\el

\bpro
As in the proof of Lemma~\ref{L:inclim} it suffices to prove $\ov\nu_{\rm Y}(\al,\de)$ is the only cluster point. By Lemma~\ref{L:liminv} each cluster point $\nu$ is an invariant law. It is also clearly homogeneous. By Lemma~\ref{L:monnu} $\nu\geq\ov\nu_{\rm Y}(\al_n,\de_n)$ and hence $\nu(\{\emptyset\})\leq\ov\nu_{\rm Y}(\al_n,\de_n)(\{\emptyset\})$. Here $\ov\nu_{\rm Y}(\al_n,\de_n)(\{\emptyset\})=0$ since $(\al_n,\de_n)\in A$ so $\nu(\{\emptyset\})=0$. By Theorem~\ref{T:ergo}, $\ov\nu_{\rm Y}(\al,\de)$ is the only homogeneous invariant law that gives zero probability to $\emptyset$ so we conclude that $\nu=\ov\nu_{\rm Y}(\al,\de)$.
\epro

\bp[General limits]
Let\label{P:lim} $A$ be the set defined in (\ref{alde}). Assume that $(\al_n,\de_n)\in A$ satisfy $\al_n\to\al$ and $\de_n\to\de$ for some $(\al,\de)\in A$. Then
\be
\ov\nu_{\rm Y}(\al_n,\de_n)\Asto{n}\ov\nu_{\rm Y}(\al,\de).
\ee
%where $\Rightarrow$ denotes weak convergence of probability measures on $\Hi(\La)$.
\ep

\bpro
As in the proof of Lemma~\ref{L:inclim} it suffices to prove $\ov\nu_{\rm Y}(\al,\de)$ is the only cluster point. By going to a subsequence, we can assume that we are in one of the following four cases: I.\ $\al_n\up\al$ and $\de_n\up\de$, II.\ $\al_n\up\al$ and $\de_n\down\de$, III.\ $\al_n\down\al$ and $\de_n\up\de$, IV.\ $\al_n\down\al$ and $\de_n\down\de$. Cases I and IV have been treated in Lemmas \ref{L:inclim} and \ref{L:declim}, respectively. In case~II we use Lemma~\ref{L:monnu} which says that in the stochastic order on $\Hi(\La)$,
\be\label{middle}
\ov\nu_{\rm Y}(\al,\de_n)\leq\ov\nu_{\rm Y}(\al_n,\de_n)\leq\ov\nu_{\rm Y}(\al_n,\de).
\ee
The left-hand side converges by Lemma~\ref{L:declim} and the right-hand side by Lemma~\ref{L:inclim}. By Lemma~\ref{L:insep} in the appendix, which is applicable by Propositions \ref{P:topol} and \ref{P:ord}, the set of monotone continuous functions on $\Hi(\La)$ is distribution determining, which allows us to conclude that in (\ref{middle}) also the expression in the middle converges. Case~III is similar.
\epro

\bpro[Proof of Theorem~\ref{T:tetcont}]
We claim that that the function
\be\label{fixco}
[0,1]\times\half\ni(\al,\de)\mapsto\tet_t(\al,\de):=\P[X_t\neq\un 0]
\ee
is continuous for each $t\geq 0$. Indeed, if $(\al_n,\de_n)\to(\al,\de)$, then by Lemma~\ref{L:rateconv} in the appendix, we can couple graphical representations $\om^n,\om$ with these rates in such a way that for the associated stochastic flows, for each $x\in S^\La_{\rm fin}$ there exists an $N<\infty$ such that $\Xb^n_{0,t}(x)=\Xb_{0,t}(x)$ for all $n\geq N$, which implies (\ref{fixco}). Since $\tet$ is the decreasing limit of the functions $\tet_t$ as $t\to\infty$, it must be upper semi-continuous.

Continuity of $\tet$ on the set $A$ from (\ref{alde}) follows from Proposition~\ref{P:lim} using (\ref{nuYX}) and Lemma~\ref{L:psico}.
\epro

\appendix

\section{Appendix}

\subsection{The stochastic order}\label{A:ord}

In this appendix we collect some general facts about the stochastic order. Throughout this appendix $\Xc$ is a compact metrisable space that is equipped with a partial order $\leq$ that is \emph{compatible with the topology} in the sense that
\be\label{compat}
\big\{(x,y)\in\Xc^2:x\leq y\big\}\mbox{ is a closed subset of }\Xc^2,
\ee
where $\Xc^2$ is equipped with the product topology.

\bl[Closedness of upset and downset]
Let\label{L:updownclosed} $\Xc$ be a compact metrisable space that is equipped with a partial order $\leq$ that is compatible with the topology. Assume that $A\subset\Xc$ is closed. Then $A^\up$ and $A^\down$ are also closed.
\el

\bpro
Since $A$ is closed and the partial order is compatible with the topology, the set $B:=\{(x,y)\in A\times\Xc:x\leq y\}$ is closed, and hence by the compactness of $\Xc$ also compact. Using the fact that the continuous image of a compact set is compact, and that $A^\up$ is the image of $B$ under the map $(x,y)\mapsto y$, we see that $A^\up$ is compact and hence closed. The same argument works for $A^\down$.
\epro

We let $B(\Xc)$ denote the space of bounded Borel measurable functions $f:\Xc\to\R$. We let $B_+(\Xc)$ denote the set of $f\in B(\Xc)$ that are \emph{monotone} in the sense that $f(x)\leq f(y)$ for all $x\leq y$. We let $\Ci(\Xc)$ denote the space of continuous functions $f:\Xc\to\R$ (which are bounded since $\Xc$ is compact) and write $\Ci_+(\Xc):=\Ci(\Xc)\cap B_+(\Xc)$. We let $\Mi_1(\Xc)$ denote the space of probability measures on $\Xc$. We will need the follwing fact.

\bl[Distribution determining property]
If\label{L:insep} $\mu,\nu\in\Mi_1(\Xc)$ satisfy $\int\mu(\di x)f(x)=\int\nu(\di x)f(x)$ for all $f\in\Ci_+(\Xc)$, then $\mu=\nu$.
\el

\bpro
Let $\Fi:=\{f\in\Ci_+(\Xc):f\geq 0\}$. By \cite[Lemma~4.37]{Swa22} it suffices to show that $\Fi$ is closed under products in the sense that $f,g\in\Fi$ imply $fg\in\Fi$, and separates points in the sense that for each $x,y\in\Xc$ with $x\neq y$, there exists an $f\in\Fi$ such that $f(x)\neq f(y)$. Closedness under products is trivial. To see that $\Fi$ separates points we observe that $x\neq y$ implies that either $\{x\}^\down\cap\{y\}^\up=\emptyset$ or $\{x\}^\up\cap\{y\}^\down=\emptyset$. By symmetry we may assume that we are in the first case. By Lemma~\ref{L:updownclosed} $\{x\}^\down$ and $\{y\}^\up$ are closed. By Theorem~I.2.1 in \cite{Nac65}, $\Xc$ is a ``normal ordered topological space'' which allows us to apply the corollary to Theorem~I.3.4 in \cite{Nac65} which tells us that there exists a monotone continuous function $f:\Xc\to[0,1]$ such that $f(x)=0$ and $f(y)=1$.
\epro

\bp[Stochastic order]
For\label{P:stochord} $\mu,\nu\in\Mi_1(\Xc)$ the following conditions are equivalent:
\begin{itemize}
\item[{\rm(i)}] It is possible to couple random variables $X,Y$ with laws $\mu,\nu$ such that $X\leq Y$ a.s.
\item[{\rm(ii)}] $\dis\int\mu(\di x)f(x)\leq\int\nu(\di x)f(x)$ for all $f\in B_+(\Xc)$.
\item[{\rm(iii)}] $\dis\int\mu(\di x)f(x)\leq\int\nu(\di x)f(x)$ for all $f\in\Ci_+(\Xc)$.
\end{itemize}
\ep

\bpro
The implications (i)$\volgt$(ii)$\volgt$(iii) are trivial and the difficult implication (iii)$\volgt$(i) is proved in \cite[Theorem~II.2.4]{Lig85}.
\epro

For $\nu,\mu\in\Mi_1(\Xc)$ we write $\mu\leq\nu$ if (i)--(iii) hold. By Lemma~\ref{L:insep} this defines a partial order on $\Mi_1(\Xc)$. We call this the \emph{stochastic order}.

%\bl[Increasing limits]
%Assume\label{L:Binc} that $\mu_n\in\Mi_1(\Xc)$ satisfy $\mu_n\leq\mu_{n+1}$ $(n\geq 0)$. Then it is possible to construct random variables $X_n,X$ such that $X_n$ has law $\mu_n$ and $X_n\up X$ a.s. Letting $\mu$ denote the law of $X$, one has
%\be\label{Binc}
%\int\mu_n(\di x)f(x)\asto{n}\int\mu(\di x)f(x)\qquad\big(f\in B(\Xc)\big).
%\ee
%An analogue statement holds for decreasing sequences.
%\el
%
%\bpro
%By Proposition~\ref{P:stochord}, for each $n\geq 1$ it is possible to couple random variables $X'_{n-1}$ and $X'_n$ with laws $\mu_{n-1}$ and $\mu_n$ such that $X_{n-1}\leq X_n$ a.s. Let $P_n(x,\,\cdot\,):=\P[X'_n\in\,\cdot\,|\,X'_{n-1}=x]$ and let $(X_n)_{n\geq 0}$ be the time-inhomogeneous Markov chain with initial law $\mu_0$ and transition kernel $P_n$ in the $n$-th step. Then $X_n$ has law $\mu_n$ and $X_n\leq X_{n+1}$ a.s.\ for all $n\geq 0$. By the compactness of $\Xc$ there exists a randon variable $X$ such that $X_n\up X$ a.s. Now (\ref{Binc}) follows from bounded pointwise convergence.
%\epro

A Feller process $(X_t)_{t\geq 0}$ with transition kernels $(P_t)_{t\geq 0}$ is called \emph{monotone} if $P_tf\in\Ci_+(\Xc)$ for all $f\in\Ci_+(\Xc)$. Equivalently, this says that $\P^x[X_t\in\,\cdot\,]\leq\P^y[X_t\in\,\cdot\,]$ for each $x\leq y$.

\bp[Upper invariant law]
Assume\label{P:upp} that $\Xc$ possesses a greatest element $\top$. Let $(P_t)_{t\geq 0}$ be the semigroup of a monotone Feller process $(X_t)_{t\geq 0}$ with state space $\Xc$. Then there exists an invariant law $\ov\nu$ of $(X_t)_{t\geq 0}$ that is uniquely characterized by the property that $\nu\leq\ov\nu$ for each invariant law $\nu$ of $(X_t)_{t\geq 0}$. Moreover, one has
\be
P_t(\top,\,\cdot\,)\Asto{t}\ov\nu,
\ee
where $\Rightarrow$ denotes weak convergence of probability measures on $\Xc$.
\ep

\bpro
This is stated for $\Xc=\{0,1\}^\Lambda$ in \cite[Theorem~III.2.3]{Lig85} and \cite[Theorem~5.4]{Swa22}, but the proof carries over to the more general setting without a change.
\epro

\subsection{Graphical representations}\label{A:partic}

In this appendix we collect some general facts about interacting particle systems and their construction from graphical representations. In particular, we provide proofs for Theorem~\ref{T:graph} and Proposition~\ref{P:perturb}. Our main reference is \cite{Swa22}.\med

\bpro[Proof of Theorem~\ref{T:graph}]
Condition (\ref{downsum}) implies condition (4.15) of \cite{Swa22}. In view of this, \cite[Thm~4.19]{Swa22} implies existence and uniqueness of solutions to (\ref{evo}), and \cite[Thm~4.20]{Swa22} implies that the process in (\ref{Xflow}) is a Feller process. By \cite[Thm~4.30]{Swa22}, the generator of this Feller process is the closure of the operator $G$ from (\ref{Gdef}), which is initially defined for functions depending on finitely many coordinates only.
\epro

We let $\Pc(\La)$ denote the set of subsets of $\La$ and write $\Pc_{\rm fin}(\La):=\{A\in\Pc(\La):|A|<\infty\}$. We equip $\Pc_{\rm fin}(\La)$, which is countable, with the discrete topology. The following lemma prepares for the proof of Proposition~\ref{P:perturb}.

\bl[Evolving set process]
Assume\label{L:xievo} (\ref{downsum})~(i) and (\ref{upsum}). Then almost surely, for each $s\in\R$ and $A\in\Pc_{\rm fin}(\La)$, there exists a unique cadlag function $\xi^{s,A}:[s,\infty)\to\Pc_{\rm fin}(\La)$ such that $\xi^{s,A}_s=A$ and
\be\label{xievo}
\xi^{s,A}_t=\left\{\ba{ll}
\big\{j\in\La:\exists i\in\xi^{s,A}_{t-}\mbox{ s.t.\ }(i,j)\in\Ri(m)\big\}
\quad&\mbox{if }(m,t)\in\om,\\[5pt]
\xi^{s,A}_{t-}\quad&\mbox{otherwise}
\ea\right.\qquad(t>s).
\ee
\el

\bpro
If $s$ and $A$ are deterministic, then $(\xi^{s,A}_{s+t})_{t\geq 0}$ is a continuous-time Markov chain with countable state space $\Pc_{\rm fin}(\La)$. Using (\ref{downsum})~(i), one can show that this process is well-defined in the sense that all its jump rates are finite, and using moreover (\ref{upsum}) one can show that this process is nonexplosive. This is extremely similar to \cite[Lemma~4.21]{Swa22}, except that the process $(\zeta_{u-t}(A))_{t\geq 0}$ considered there ``runs backwards in time'' and (\ref{upsum}) is replaced by (\ref{downsum})~(ii). Since the proof is essentially the same, we omit the details. Using an argument as in \cite[Lemma~4.22]{Swa22}, one can remove the assumption that $s$ and $A$ are deterministic, i.e., one can show that the statement holds a.s.\ for all $s\in\R$ and $A\in\Pc_{\rm fin}(\La)$ simultaneously.
\epro

\bpro[Proof of Proposition~\ref{P:perturb}]
Fix $s\in\R$ and $x,y\in S^\La$ such that $A:=\{i\in\La:x(i)\neq x(j)\}$ is finite, and set $X^{s,x}_t:=\Xb_{s,t}(x)$ $(t\geq 0)$. Conditions (\ref{downsum})~(i) and (\ref{upsum}) allow us to apply Lemma~\ref{L:xievo} which says that (\ref{xievo}) has a unique solution $(\xi^{s,A}_t)_{t\geq 0}$. Let
\be\ba{r@{\,}l}\label{omprim}
\dis\om':=&\dis\big\{(m,t)\in\om:t>s,\ \Di(m)\cap\xi^{s,A}_{t-}\neq\emptyset\big\}\\[5pt]
&\dis\cup\,\big\{(m,t)\in\om:t>s,\ \exists(i,j)\in\Ri(m)\mbox{ s.t.\ }i\in\xi^{s,A}_{t-},\ j\in\Di(m)\big\}.
\ec
%\bc\label{omprim}
%\dis\om'&:=&\dis\big\{(m,t)\in\om:t>s,\ \Di(m)\cap\xi^{s,A}_{t-}\neq\emptyset\big\}\\[5pt]
%&&\dis\cup\,\big\{(m,t)\in\om:t>s,\ \exists(i,j)\in\Ri(m)\mbox{ s.t.\ }i\in\xi^{s,A}_{t-},\ j\in\Di(m)\big\}.
%\ec
It follows from (\ref{downsum})~(i), (\ref{upsum}), and the finiteness of $(\xi^{s,A}_t)_{t\geq 0}$ that we can order the elements of $\om'$ as
\be
\om'=\big\{(m_k,t_k):k\geq 1\big\}\quad\mbox{with}\quad t_1<t_2<\cdots.
\ee
We now define $(X^{s,y}_t)_{t\geq s}$ by first setting
\be
X^{s,y}_t(i):=X^{s,x}_t(i)\quad(t\geq s,\ i\not\in\xi^{s,A}_t),
\ee
and then setting
\be\ba{r@{\,}c@{\,}ll}
\dis X^{s,y}_t(i)&:=&\dis y(i)\quad&\dis(t\in[s,t_1),\ i\in A),\\[5pt]
\dis X^{s,y}_t(i)&:=&\dis m_k(X^{s,y}_{t_k-})(i)\quad&\dis(t\in[t_k,t_{k+1}),\ i\in \xi^{s,A}_{t_k},\ k\geq 1).
\ec
It is then straightforward to check that $(X^{s,y}_t)_{t\geq s}$ solves (\ref{evo}) with initial condition $X^{s,y}_s=y$. Using also condition (\ref{downsum})~(ii) (which has not been used up to this point) we can apply Theorem~\ref{T:graph} to conclude that (\ref{evo}) has a unique solution and $X^{s,y}_t=\Xb_{s,t}(y)$. Then
\be
\big\{i\in\La:\Xb_{s,t}(x)(i)\neq\Xb_{s,t}(y)(i)\big\}\sub\xi^{s,A}_t\qquad(t\geq s),
\ee
so the sets in (\ref{difset}) are finite for all $t\geq s$.
\epro

We need some approximation results that allow us to conclude that particle systems can be approximated by finite systems, or by systems whose rates converge in an appropriate way. In the setting of Theorem~\ref{T:graph}, if $\om'$ is a finite subset of $\Gi\times\R$ such that no two elements of $\om'$ have the same time coordinate, then we can define a stochastic flow $(\Xb^{\om'}_{s,t})_{s\leq t}$ by setting, for each $s\leq u$,
\be\ba{c}
\dis\Xb^{\om'}_{s,u}:=m_n\circ\cdots\circ m_1\quad\mbox{where}\quad\\[5pt]
\dis\big\{(m,t)\in\om':s<t\leq u\big\}=\big\{(m_1,t_1),\ldots,(m_n,t_n)\big\}
\quad\mbox{with}\quad t_1<\cdots<t_n.
\ec
We cite the following result from \cite[Lemma~4.24]{Swa22}.

\bp[Finite approximation]
Under\label{P:finap} the assumptions of Theorem~\ref{T:graph}, almost surely for all $s\leq u$ and for each sequence $\om_n$ of finite subsets of $\om$ such that $\om_n\up\om$, one has
\be
\Xb^{\om_n}_{s,t}(x)\asto{n}\Xb_{s,t}(x)\qquad(x\in S^\La),
\ee
where $\to$ denotes convergence in the product topology.
\ep

Let $\La$ be countable, let $S$ be finite, and let $\Gi$ be a countable collection of local maps $m:S^\La\to S^\La$. Assume (\ref{mas}). Let $\Rb$ be the set of all collections $\rb=(r_m)_{m\in\Gi}$ of nonnegative rates that satisfy (\ref{downsum}) and (\ref{upsum}). By Lemma~\ref{L:finsys}, for each $\rb\in\Rb$, we can construct a stochastic flow $(\Xb_{s,t})_{s\leq t}$ that maps the space $S^\La_{\rm fin}$ into itself. We need a result that says that if a sequence of rates $\rb^n$ converges in an appropriate sense, then the associated stochastic flows converge.

\bl[Convergence of finite systems]
For\label{L:rateconv} each $i\in\La$, set
\be
\Gi_i:=\big\{m\in\Gi:i\in\Di(m)\big\}\cup\big\{m\in\Gi:\exists j\in\Di(m)\mbox{ s.t.\ }(i,j)\in\Ri(m)\big\}.
\ee
Assume that $\rb^n,\rb\in\Rb$ satisfy
\be\label{rateconv}
\sum_{m\in\Gi_i}|r^n_m-r_m|\asto{n}0\quad\forall i\in\La.
\ee
Then it is possible to couple the graphical representations $\om^n,\om$ with rates $\rb^n,\rb$ in such a way that the associated stochastic flows satisfy
\be
\Xb^n_{s,u}(x)\asto{n}\Xb_{s,u}(x)\quad{\rm a.s.}\quad\forall s\leq u,\ x\in S^\La_{\rm fin},
\ee
where the convergence is with respect to the discrete topology on $S^\La_{\rm fin}$.
\el

\bpro
Let $\pi$ be a Poisson point set on $\Gi\times\R\times\half$ with intensity measure
\be
\mu\big(\{m\}\times[s,t]\times[0,r]\big):=r(t-s)\qquad(m\in\Gi,\ s\leq t,\ r\geq 0).
\ee
Then for each collection of rates $\rb=(r_m)_{m\in\Gi}$ we can define a graphical representation with these rates by setting
\be
\om:=\big\{(m,t):(m,t,r)\in\pi,\ r\leq r_m\big\}.
\ee
Let $\om^n$ and $\om$ be constructed in this way for the rates $\rb^n$ and $\rb$. Fix $s\leq u$ and $x\in S^\La_{\rm fin}$, set $\xi_t:=\{i\in\La:\Xb_{s,t}(x)(i)\neq 0\}$ $(t\geq s)$, and consider the set (compare (\ref{omprim}))
\be\ba{r@{\,}l}\label{omxi}
\dis\ti\om:=&\dis\big\{(m,t)\in\om:s<t\leq u,\ \Di(m)\cap\xi_{t-}\neq\emptyset\big\}\\[5pt]
&\dis\cup\,\big\{(m,t)\in\om:s<t\leq u,\ \exists(i,j)\in\Ri(m)\mbox{ s.t.\ }i\in\xi_{t-},\ j\in\Di(m)\big\}.
\ec
Define $\ti\om^n$ similarly, with $\om$ in (\ref{omxi}) replaced by $\om^n$ (but still using the same process $(\xi_t)_{t\geq 0}$ which is defined in terms of $\om$). The condition (\ref{rateconv}) guarantees that almost surely $\ti\om^n=\om^n$ for all $n$ large enough, and hence also $\Xb^n_{s,u}(x)=\Xb_{s,u}(x)$ for all $n$ large enough.
\epro

\end{document}